\documentclass{amsart}

\usepackage{amsfonts,amssymb,verbatim,amsmath,amsthm,latexsym,textcomp,amscd}
\usepackage{latexsym,amsfonts,amssymb,epsfig,verbatim}
\usepackage{amsmath,amsthm,amssymb,latexsym,graphics,textcomp}
\usepackage{graphicx}
\usepackage{color}
\usepackage{url}
\usepackage{stmaryrd}
\usepackage{enumitem}
\usepackage{tikz}
\usepackage{mathtools}
\usetikzlibrary{positioning}
\usepackage{hyperref}
\hypersetup{
	colorlinks   = true, urlcolor     = blue, linkcolor    = black, citecolor    = black}
\begin{document}

	\newtheorem{theorem}{Theorem}[section]
	\newtheorem{proposition}[theorem]{Proposition}
	\newtheorem{lemma}[theorem]{Lemma}
	\newtheorem{corollary}[theorem]{Corollary}
	\newtheorem{definition}[theorem]{Definition}
	\newtheorem{notation}[theorem]{Notation}
	\newtheorem{convention}[theorem]{Convention}
	\newtheorem{conjecture}[theorem]{Conjecture}
	\newtheorem{claim}[theorem]{Claim}
	\newtheorem{question}[theorem]{Question}
	\newtheorem{remark}{Remark}
	\newtheorem{example}{Example}
	\newcommand{\map}{\rightarrow}
	\newcommand{\C}{\mathbb C}
	\newcommand\AAA{{\mathcal A}}
	\def\AA{\mathcal A}
	
	\def\L{{\mathcal L}}
	\def\al{\alpha}
	\def\A{{\mathcal A}}
	\newcommand\BB{{\mathcal B}}
	\newcommand\CC{{\mathcal C}}
	\newcommand\DD{{\mathcal D}}
	\newcommand\EE{{\mathcal E}}
	\newcommand\FF{{\mathbb F}}
	\newcommand\GG{{\mathcal G}}
	\newcommand\GB{{\mathbb G}}
	\newcommand\HH{{\mathcal H}}
	\newcommand\II{{\mathcal I}}
	\newcommand\JJ{{\mathcal J}}
	\newcommand\KK{{\mathcal K}}
	\newcommand\LL{{\mathbb L}}
	\newcommand\LS{{\mathcal L}} %S for script
	\newcommand\MM{{\mathcal M}}
	\newcommand\NN{{\mathbb N}}
	\newcommand\OO{{\mathcal O}}
	\newcommand\PP{{\mathcal P}}
	\newcommand\QQQ{{\mathbb Q}}
	\newcommand\QQ{{\mathcal Q}}
	\newcommand\RR{{\mathbb R}}
	\newcommand\SSS{{\Sigma}}
	\newcommand\TT{{\mathcal T}}
	\newcommand\UU{{\mathcal U}}
	\newcommand\VV{{\mathcal V}}
	\newcommand\WW{{\mathcal W}}
	\newcommand\XX{{\mathcal X}}
	\newcommand\YY{{\mathcal Y}}
	\newcommand\ZZ{{\mathbb Z}}
	\newcommand\hhat{\widehat}
	\newcommand{\Lam}[1]{\ensuremath{\partial^{(2)}_{#1}}}
	\newcommand\XH{\mathcal{X}}
	\newcommand{\YH}{\mathcal{Y}}
	\newcommand{\mh}{\hat{d}_1}
	\newcommand{\X}{\widehat{X}}
	\newcommand{\Y}{\widehat{Y}}
	\newcommand{\F}{\widehat{F}}
	\newcommand{\m}{\hat{d}}
	\newcommand{\CV}{\text{CV}}
	\newcommand{\stab}{\text{Stab}}
	\newcommand{\f}{\hat{f}}
	\newcommand{\T}{\widehat{T}}
	\newcommand{\hsigma}{\hat{\sigma}}
	\newcommand{\hgamma}{\hat{\gamma}}
	\newcommand{\hrho}{\hat{\rho}}
	\newcommand{\secref}[1]{Section~\ref{#1}}
	\newcommand{\thmref}[1]{Theorem~\ref{#1}}
	\newcommand{\lemref}[1]{Lemma~\ref{#1}}
	\newcommand{\rmkref}[1]{Remark~\ref{#1}}
	\newcommand{\propref}[1]{Proposition~\ref{#1}}
	\newcommand{\corref}[1]{Corollary~\ref{#1}}
	\newcommand{\probref}[1]{Problem~\ref{#1}}
	\newcommand{\eqnref}[1]{~{\textrm(\ref{#1})}}
	\newcommand{\boldG}{\mbox{{\bf G}}}
	
	\def\Ga{\Gamma}
	\def\Z{\mathbb Z}
	
	\def\diam{\operatorname{diam}}
	\def\dist{\operatorname{dist}}
	\def\hull{\operatorname{Hull}}
	
	\def\length{\operatorname{length}}
	\newcommand\RED{\textcolor{red}}
	\newcommand\GREEN{\textcolor{green}}
	\newcommand\BLUE{\textcolor{blue}}
	\def\mini{\scriptsize}
	
	\def\acts{\curvearrowright}
	\def\embed{\hookrightarrow}
	
	\def\ga{\gamma}
	\newcommand\la{\lambda}
	\newcommand\eps{\epsilon}
	\def\geo{\partial_{\infty}}
	\def\bhb{\bigskip\hrule\bigskip}
	
	\title[]{Relative hyperbolicity of ascending HNN extension of groups}
	
	\author{Swathi Krishna}
	\address{Kerala School of Mathematics, Kozhikode, Kerala, India}
	\email{swathi280491@gmail.com}
	\date{\today}
	
	\begin{abstract}
		We prove that for a finitely generated group $G$ with a free factor system $\GG$ and an injective endomorphism $\phi$ that preserves $\GG$, then the ascending HNN extension $G\ast_{\phi}$ is hyperbolic relative to a collection of maximal parabolic subgroups. As a corollary, we see if the injective endomorphism, $\phi$, of a finite rank free group $\FF$ is exponentially growing, the ascending HNN extension $\FF\ast_{\phi}$ is relatively hyperbolic.
	\end{abstract}
	
	\maketitle
	
	%\tableofcontents
	\section{Introduction}
	\noindent
	Let $\FF$ be a finite rank free group and $\phi:\FF\to\FF$ be an automorphism. The geometry of $\FF\rtimes_{\phi}\Z$ has been studied extensively, mainly using the theory of train tracks developed by Bestvina and Handel in \cite{traintrack}. 
	In \cite{brinkmann}, Brinkmann showed that $\FF\rtimes_{\phi}\mathbb{Z}$ is hyperbolic if and only if $\phi$ is atoroidal. Ghosh generalized this in \cite{pritam}, where he showed that $\FF\rtimes_{\phi}\mathbb{Z}$ is relatively hyperbolic if and only if $\phi$ is exponentially growing. Much more generally, in \cite{DahmaniLi}, Dahmani and Li proved a combination theorem for relative hyperbolcity of free product-by-cyclic groups, followed by the work of Dahmani and Suraj Krishna (see \cite{surajdahmani}) for hyperbolic-by-cyclic groups. \\
	Let $G$ be a finitely generated group and $\phi:G\to G$ be an injective endomorphism. The mapping torus of $\phi$ is an ascending HNN extension $$G\ast_{\phi}=\langle G,t\mid tgt^{-1}=\phi(g), g\in G\rangle.$$
	
	\noindent
	In \cite{IKendo}, Ilya Kapovich showed that, when $\phi:\FF\to\FF$ is an immersion, $\FF\ast_{\phi}$ is hyperbolic if and only if $\FF\ast_{\phi}$ does not contain the Baumslag-Solitar group $BS(1,d)=\langle a,t\mid tat^{-1}=a^d\rangle$, $d\geq1$. This is equivalent to saying that there exists no $k,d>0$, $g\in\FF\setminus\{1\}$ such that $g^d$ and $\phi^k(g)$ are conjugate. In \cite{JPMhyperbolic,JPMrel}, Mutanguha generalized this for any injective endomorphism of $\FF$. Using the techniques developed in \cite{JPMrel}, he proved the relative hyperbolicity of $\FF\ast_{\phi}$ in \cite{JPMrelhyp}. 
	
	Now, let $G=H_1\ast H_2\ast\cdots\ast H_p\ast\FF_k$, $\GG=\{[H_1],\ldots,[H_p]\}$, with Scott complexity $(k,p)\neq(0,1)$, and let $\phi:(G,\GG)\to(G,\GG)$ be an injective endomorphism. The following is our main result.
	
	\vspace{2mm}
	\noindent
	{\bf Theorem \ref{main-poly}.}
	Let $G$ be a finitely generated group with a free factor system $\GG$ with Scott complexity $(k,p)\neq(0,1)$, and $\phi:(G,\GG)\to(G,\GG)$ be an injective, non-surjective endomorphism. Let $\PP$ be the collection of conjugacy classes of maximal polynomially growing subgroups for $\phi$ on a (any) tree in $\TT_{\GG}$. If $\phi$ is strictly type preserving relative to $\PP$, then the ascending HNN extension $G\ast_{\phi}=\langle G,t\mid t^{-1}gt=\phi(g),g\in G\rangle$ is hyperbolic relative to a collection of ascending HNN extension of some elements in $\PP$.
	
	\vspace{2mm}
	\noindent
	As a corollary, we get the following.
	
	\vspace{2mm}
	\noindent
	{\bf Theorem \ref{main-free}.}
	Let $\phi:\FF\to\FF$ be an injective non-surjective endomorphism of a finite rank free group. Let $\phi$ be exponentially growing and $\PP$ be the collection of conjugacy classes of maximal polynomially growing subgroups for $\phi$. If $\phi$ is strictly type preserving relative to $\PP$, then $\FF\ast_{\phi}=\langle\FF,t\mid t^{-1}gt=\phi(g),g\in\FF\rangle$ is hyperbolic relative to a collection of ascending HNN extension of some elements in $\PP$.
	
	\vspace{2mm}
	\noindent
	As a consequence of the results known for $\FF\ast_{\phi}$, we give the necessary and sufficient condition for $\FF\ast_{\phi}$ to have the rapid decay property.
	
	\vspace{2mm}
	\noindent
	{\bf Theorem \ref{rapiddecay}.}
	Let $\phi:\FF\to\FF$ be an injective endomorphism. Then $\FF\ast_{\phi}$ does not have the rapid decay property if and only if $\FF\ast_{\phi}$ contains $BS(1,d)$ for some $d\geq 2$.
	%\vspace{2mm}
	%\noindent
	%{\bf Outline of the paper:} In Section \ref{section 2}, we recall some definitions and results related to relative hyperbolicity and tree of spaces. In Section \ref{section 3}, 
	
	%\medskip
	%\noindent 
	%{\bf Acknowledgments.} 
	%The author is grateful to for pointing out a mistake and for the helpful comments.
	
	\section{Preliminaries}\label{section 2}
	\noindent
Let $G=H_1\ast H_2\ast\cdots\ast H_p\ast\FF_k$, and $\GG=\{[H_1],\ldots,[H_p]\}$ be the free factor system.
\begin{notation}
	\begin{enumerate}
		\item An injective endomorphism $\phi:(G,\GG)\to(G,\GG)$ is an injective endomorphism of $G$ such that for every $1\leq i\leq p$, there exists $1\leq j\leq p$ and $x\in G$ such that $\phi(H_i)\leq x\phi(H_j)x^{-1}$. 
		\item The outer class of $\phi$ is denoted by $[\phi]$.
		\item For any element $g$ (resp. free factor $H$) in $G$, $[g]$ (resp. $[H]$) denotes its conjugacy class in $G$.
		\item For a pair of collection of subgroups $\AA$ and $\BB$ of $G$, $\AA\preccurlyeq\BB$ if for every $A\in\AA$, there exists $B\in\BB$, $x\in G$ such that $A\leq xBx^{-1}$.
		\item $\AA$ is $[\phi]$-invariant if $\phi(\AA)\preccurlyeq\AA$.
		Further, $\AA$ is $[\phi]$-fixed if $\phi(\AA)=\AA$.
		\item $Hyp(\GG)$ denotes the set of all $g\in G$ such that there exists no $[H_i]\in\GG$ and $x\in G$ with $g\in xH_ix^{-1}$.
		\item Let $G$ act on a tree $T$ by isometries. For a finitely generated subgroup $H\leq G$ containing a hyperbolic isometry, $T_H\subset T$ denotes the minimal $H$-invariant subtree. In fact, $T_H$ is the union of axes of all hyperbolic elements of $H$ in $T$.
	\end{enumerate}   
\end{notation}
%\begin{definition}
%For $i=1,2$, let $\HH_i=\{[H^i_1],\ldots,[H^i_{k_i}]\}$ be free factor systems of $G$ such that $G=H^i_1\ast\cdots\ast H^i_k\ast\FF_{r_i}$. 
%$\HH_1$ is {\bf lower than} $\HH_2$ if for any $[H^1_j]\in\HH_1$, there exists $[H^2_k]\in\HH_2$ and $g\in\FF$ such that $H^1_j\leq gH^2_kg^{-1}$. We denote this by $\HH_1\preccurlyeq\HH_2$. If $\HH_1$ is strictly lower than $\HH_2$, we denote it by $\HH_1\prec\HH_2$.
%\end{definition}
\begin{lemma}\cite[Lemma 1.1]{DahmaniLi}\label{scott_complexity}
	If $\HH\prec\HH'$, then the Scott complexity of $\HH'$ is strictly smaller than that of $\HH$ in the lexicographic order. For any $[H]\in\HH'$, the Scott complexity of the free factor system induced by $\HH$ on $H$ is also strictly smaller than that of $\HH'$ in the lexicographic order.
\end{lemma}
%\begin{proof}
%    Let $\HH=\{[H_1],\ldots,[H_k]\}, \HH'=\{[H'_1],\ldots,[H'_l]\}$. Then $G=H_1\ast\cdots\ast H_k\ast\FF_p$ and $G=H'_1\ast\cdots\ast H'_l\ast\FF_q$. Since $\HH\prec\HH'$, for each $1\leq i\leq l$, there exists $i_1,\ldots,i_{k_i}\in\{1,\lodts,k\}$ such that $H'_i=H_{i_1}\ast\cdots\ast H_{i_{k_i}}\ast\FF_{q_i}$. Then clearly, $k\geq l$ and $q_i\ast q$ for each $i$. Therefore, $(p,k)$ is smaller than $(q,l)$.
%Also, for each $i$, $(q_i,i_{k_i})$ is smaller than $(q,l)$.
%\end{proof}

\begin{definition}\label{defn-irred1}\begin{enumerate}
		\item $\phi:(G,\GG)\to(G,\GG)$ is \textbf{atoroidal relative to $\GG$} if there exists no $g\in Hyp(\GG)$, $k,d\geq 1$ such that $\phi^k([g])=[g^d]$.
		\item $\phi$ is \textbf{irreducible relative to $\GG$} if there exists no proper free factor system $\GG'\succ\GG$ that is invariant under $\phi$.
		\item $\phi$ is \textbf{fully irreducible relative to $\GG$} if for each $i>0$, $\phi^i$ is irreducible relative to $\GG$. 
	\end{enumerate}
\end{definition}
\begin{remark}
	An injective endomorphism $\phi$ is atoroidal if there exists no $g\in G\setminus\{1\}$, $k,d\geq 1$ such that $\phi^k([g])=[g^d]$. Note that when $\phi$ is an automorphism, the above definition of $\phi$ being atoroidal is equivalent to the original definition of atoroidality, i.e., there exists no $g\in G\setminus\{1\}$, $k\geq 1$ such that $\phi^k(g)$ is conjugate to $g$.
\end{remark}
\subsection{Relative hyperbolicity}
Here we refer the reader to \cite{farb}, \cite{bowditch-relhyp}, \cite{groves-manning} and \cite{osin} for various equivalent definitions and results of relatively hyperbolic groups. We recall the following construction from \cite{farb}. Let $X$ be a graph and $\mathcal{Y}=\{Y_\alpha\}_{\alpha\in\Lambda}$ be a collection of uniformly separated disjoint subgraphs. Corresponding to each $Y_{\alpha}$, we attach a new vertex $\nu(Y_\alpha)$ and join each vertex in $Y_\alpha$ to $\nu(Y_\alpha)$ by an edge of length $\frac{1}{2}$. These new vertices are called the {\em cone points} and this construction is called the {\bf coning-off} of $X$. This new coned-off graph is denoted by $\mathcal{E}(X,\mathcal{Y})$ or $\widehat{X}$.
\begin{remark}
	Let $G$ be a finitely generated group and $H$ be a finitely generated subgroup. Let $X$ denote the Cayley graph of $G$ with respect to a finite generating set and $\mathcal{Y}$ denote the collection of all subgraphs corresponding to the cosets of $H$ in $G$. Then by the definition of Farb, $G$ is weakly hyperbolic relative to $H$ if $\mathcal{E}(X,\mathcal{Y})$ is a hyperbolic metric space.
\end{remark}
\noindent
Throughout the paper, we refer to the strong relative hyperbolicity simply as relative hyperbolicity.
\subsection{Tree of spaces}
\begin{definition}
	A {\bf graph of groups} $(\mathcal{G},Y)$ consists of a finite graph $Y$ with vertex set $V$ and edge set $E$ and for each vertex $v \in V$, there is a group $G_v$ (vertex group) and for each edge $e \in E$, there is a group $G_e$ (edge group), along with the monomorphisms: 	
	$$\phi_{o(e)}:G_e \to G_{o(e)},\phi_{t(e)}:G_e \to G_{t(e)},$$	
	with the extra condition that $G_{\bar{e}}=G_e$.   
\end{definition}
%\begin{definition}
%A graph of groups $(\mathcal{G},Y)$ is a {\bf graph of relatively hyperbolic groups} if for each $v \in V(Y)$, $G_v$ is hyperbolic relative to a collection of subgroups $\{H_{v,\alpha}\}_{\alpha}$ and for each $e \in E(Y)$, $G_e$ is hyperbolic relative to a collection of subgroups $\{H_{e,\alpha}\}_{\alpha}$.   
%\end{definition}
\noindent
Let $G$ be a finitely generated group with a graph of groups decomposition $(\mathcal{G},Y)$. Let $\TT$ be its Bass Serre tree. Recall that $V(\TT)=\bigsqcup_{v\in V(Y)}G/G_v$ and $E(\TT)=\bigsqcup_{e\in E(Y)}G/\phi_{t(e)}(G_e)$. For an edge $g\phi_{t(e)}(G_e)$, $$o(g\phi_{t(e)}(G_e))=gG_{o(e)},t(g\phi_{t(e)}(G_e))=geG_{t(e)}.$$ One can further associate a tree of spaces with it. For each $e\in E(Y)$, we fix a finite generating $S_e$ of $G_e$ and for each $v\in V(Y)$, we fix a finite generating set $S_v$ of $G_v$ such that for an edge $e$ with $o(e)=v$, $\phi_{o(e)}(S_e)\subset S_v$. Let $X_v$ (resp. $X_e$) denote the Cayley graph of $G_v$ (resp. $G_e$) with respect to $S_v$ (resp. $S_e$). Let $S=\bigcup_{v\in V(Y)}S_v\bigcup\big(E(Y)\setminus E(T_Y)\big)$, where $T_Y$ is a maximal subtree in $Y$.
\begin{definition}\cite{mit0}
	A {\bf tree of spaces} is a metric space $X$ admitting the map $p:X\to \TT$ such that: 
	\begin{enumerate}
		\item For $\tilde{v}=gG_v\in V(\TT)$, $X_{\tilde{v}} = p^{-1}(\tilde{v})$ is a subgraph of $\Gamma(G,S)$ with $V(X_{\tilde{v}}) = gG_v$ and $gx, gy \in X_{\tilde{v}}$ are connected by an edge if $x^{-1}y \in S_v$. 	
		\item For $\tilde{e}=g\phi_{t(e)}(G_e)\in E(\TT)$, $X_{\tilde{e}} = p^{-1}(\tilde{e})$ is a subgraph of $\Gamma(G,S)$ with $V(X_{\tilde{e}}) = ge\phi_{t(e)}(G_e)$ and $gex, gey \in X_{\tilde{e}}$ are connected by an edge if $x^{-1}y \in \phi_{t(e)}(S_e)$. 		
		\item For an edge $\tilde{e} = g\phi_{t(e)}(G_e)$ connecting vertices $\tilde{u} = gG_{o(e)}$ and $\tilde{v} = geG_{t(e)}$, if $x \in\phi_{t(e)}(G_e)$, we join $gex \in X_{\tilde{e}}$ to $gexe^{-1} \in X_{\tilde{u}}$ and $gex \in X_{\tilde{v}}$ by edges of length $\frac{1}{2}$. 
		These extra edges induce maps $f_{\tilde{e},\tilde{u}} : X_{\tilde{e}} \to X_{\tilde{u}}$ and $f_{\tilde{e},\tilde{v}} : X_{\tilde{e}} \to X_{\tilde{v}}$ with $f_{\tilde{e},\tilde{u}}(gex) = gexe^{-1}$ and $f_{\tilde{e},\tilde{v}}(gex) = gex$.
\end{enumerate}\end{definition}
\begin{definition}\cite{MjReeves}
	$(\GG,Y)$ is a {\bf graph of relatively hyperbolic groups} if for each $v\in V(Y)$, there exists a finite collection of subgroups $\HH_v=\{H_{v,\alpha}\}_{\alpha\in\Lambda}$ of $G_v$ such that $G_v$ is hyperbolic relative to $\HH_v$, and similarly, for each $e\in E(Y)$, there exists a finite collection of subgroups $\HH_e=\{H_{e,\alpha}\}_{\alpha\in\Lambda}$ of $G_e$ such that $G_e$ is hyperbolic relative to $\HH_e$.     
\end{definition}
\begin{definition}\cite{MjReeves}\label{conditions}
	\begin{enumerate}
		\item \textbf{QI-embedded condition:} A graph of groups $(\mathcal{G},Y)$ is said to satisfy the {\em qi (quasiisometric)-embedded condition} if for every $e \in E(Y)$, the monomorphisms $\phi_{o(e)}$ and $\phi_{t(e)}$ are qi-embeddings.\\
		This is equivalent to saying that for each $\tilde{e}\in E(\TT)$, the maps $f_{\tilde{e},o(\tilde{e})}$ and $f_{\tilde{e},t(\tilde{e})}$ are qi embedding, and so the tree of spaces $X$ also satisfies a qi-embedded condition.
		\item \textbf{Strictly type-preserving:} A graph of relatively hyperbolic groups is {\em strictly type-preserving} if for every $e \in E(Y)$, each $\phi_{o(e)}^{-1}(H_{v,\alpha})$ and $\phi_{t(e)}^{-1}(H_{v,\alpha})$ is either empty or some $H_{e,\alpha}$.
		\item \textbf{QI-preserving electrocution condition:} $(\mathcal{G},Y)$ satisfies {\em qi-preserving electrocution condition} if induced maps $\hat{\phi}_{o(e)}: \widehat{\Gamma}_{e} \to \widehat{\Gamma}_{o(e)}$ and $\hat{\phi}_{t(e)}: \widehat{\Gamma}_{e} \to \widehat{\Gamma}_{t(e)}$ are uniform qi-embeddings. Here, $\widehat \Gamma_{e}$, $\widehat \Gamma_{o(e)}$ and $\widehat \Gamma_{t(e)}$ denote the coned-off Cayley graphs of $G_e$, $G_{o(e)}$ and $G_{t(e)}$ respectively relative to the corresponding maximal parabolic subgroups. 
\end{enumerate}\end{definition}
\noindent
Given a graph of relatively hyperbolic groups $(\mathcal{G},Y)$ that satisfies the strictly type-preserving condition, we can cone-off the vertex spaces and edge spaces in $X$ and join the cone points by edges according to the maps induced by $f_{e,v}$ to get an induced tree of coned-off spaces. We denote this by $\mathcal{TC}(X)$. If the qi-embedded and qi-preserving electrocution conditions are also satisfied, $\mathcal{TC}(X)$ is a tree of spaces satisfying a qi-embedded condition.\\
The subgraph of $\mathcal{TC}(X)$ with the set of cone points as the vertex set and corresponding edges is called the {\bf cone locus}. This is in fact a forest and each connected component can be identified with a subtree in $\mathcal{T}$. We denote the collection of maximal connected components of the cone locus by $\YY=\{Y_{\alpha}\}_{\alpha\in\Lambda}$. Each $Y_{\alpha}$ gives rise to a tree of horophere-like sets over it. We denote this collection by $\CC=\{C_\alpha\}_{\alpha\in\Lambda}$. One is referred to \cite[Subsection 1.4]{MjPal} for more details.

	\section{Setup}\label{section 3}
	\noindent
{\bf Notation:} Given a path $\rho$ in a tree $T$, we denote the tightened path or the geodesic joining the endpoints of $\rho$ by $[\rho]$.

\vspace{2mm}
\noindent
Let $G=H_1\ast\cdots\ast H_p\ast\FF_k$ and $\GG=\{[H_1],\ldots,[H_p]\}$. 
Let $\phi:(G,\GG)\to(G,\GG)$ be an injective endomorphism that is strictly type-preserving.
Recall that $T$ is a $G$-tree if $G$ acts on $T$ by isometries. Following \cite{gl-outer}, let $\TT'_{\GG}$ be a collection of metric simplicial $G$-trees such that:\\
(0) $T$ has no redundant vertices ($v\in V(T)$ is redundant if it has valence $2$ and if $g\in G$ fixes $v$, then it fixes every edge incident on $v$).\\
(1) Action is non-trivial and minimal.\\
(2) Edge stabilizers are trivial and there are finitely many edge orbits,\\
(3) For $1\leq i\leq p$, for each conjugate of $H_i$ in $G$, there exists a unique vertex fixed by this conjugate. Other vertices have trivial stabilizer and are called {\em free vertices}.\\
\noindent
Then the {\bf outer space} $\TT_{\GG}$ of $(G,\GG)$ is a collection of equivalence classes of trees in $\TT'_{\GG}$, where $T\sim T'$ if there exists a homothety $f:T\to T'$.
%Both $\TT'_{\GG}$ and $\TT_{\GG}$ have length function topology.
\begin{remark}
	Let $T\in\TT_{\GG}$ and let $\alpha:G\to\text{Isom}(T)$ denote the action of $G$ on $T$. When $\phi$ is an automorphism, $(T,\alpha)\cdot\phi=(T,\alpha\circ\phi)$. But when $\phi$ is an injective, non-surjective endomorphism, $(T,\alpha)\cdot\phi=(T_{\phi(G)},\alpha\circ\phi)$. 
	%In this case, vertex stabilizers are either trivial or subgroups of conjugates of elements of $\GG$. 
\end{remark}
\subsection{Existence of train track map}
We follow the results mainly from \cite{bf} and \cite{jpmutanguha} to prove that an injective endomorphism $\phi:(G,\GG)\to(G,\GG)$ that is irreducible relative to $\GG$ has a train track representative. 
We first recall the definition of train track maps (see \cite{DahmaniLi,FranMartino}).
\begin{definition}
	\begin{itemize}
		\item Given any $v\in V(T)$, a {\em turn} at $v$ is a pair of edges $(e_1,e_2)$ such that $i(e_1)=v=i(e_2)$. 
		\item Given $T\in\TT_{\GG}$, a {\em train track structure} on $T$ is a $G$-invariant equivalence relation on the set of oriented edges at each vertex of $T$ such that there exists at least two equivalence classes at each vertex. Each equivalence class is called a {\em gate}.
		\item A turn $(e_1,e_2)$ is {\em legal} if $e_1,e_2$ are in distinct equivalence classes. Otherwise it is {\em illegal}. A path is legal if every turn in the path is legal.
		\item For $T\in\TT_{\GG}$, a $\phi$-equivariant map $f:T\to T$ is a \textbf{train track map} if $(1)$ $f$ sends edges to legal paths, and $(2)$ for any $v\in V(T)$, if $f(v)\in V(T)$, a legal turn at $v$ maps to a legal turn at $f(v)$.        
	\end{itemize}
\end{definition}
\begin{lemma}\label{top_rep}
	Given an injective non-surjective endomorphism $\phi:(G,\GG)\to(G,\GG)$ and $T\in\TT_{\GG}$, there exists a $\phi$-equivariant map $f:T\to T$ representing $\phi$. 
\end{lemma}
\noindent
Proof of this result is the same as that of \cite[Lemma 4.2]{FranMartino}, but note that if $\phi$ is not surjective, then $f$ is also not surjective. It also follows from their result that if $f'$ is another $\phi$-equivariant map representing $\phi$, then $f$ and $f'$ agree on the non-free vertices. Following \cite{jpmutanguha}, we will refer to this map as a {\bf $\GG$-relative weak representative of $\phi$}.
%\begin{definition}
%   Let $f:T\to T$, where $f,T$ are as in Lemma \ref{top_rep}.
%  \begin{enumerate}
%     \item $f$ is a {\bf $\GG$-relative representative of $\phi$} if $f$ is $\GG$-relative weak representative with no pre-trivial edges and $T$ has no bivalent vertices.
%        \item $f$ is {\bf minimal} if it is a $\GG$-relative weak representative with no orbit-closed $f$-invariant sub-forests whose components are bounded. 
%   \end{enumerate}
%\end{definition}
\noindent
Recall that a $\GG$-relative weak representative induces a gate structure on $T$, i.e., turn $(e,e')$ is illegal if $f(e)$ and $f(e')$ have the same initial segment. This is a $G$-invariant equivalence relation. As in \cite{DahmaniLi}, this is the gate structure used in this paper. 

\noindent
Let $\phi:(G,\GG)\to(G,\GG)$ be an injective endomorphism and let $f:T\to T$ be a $\GG$-relative weak representative. Let $\{e_1,\ldots,e_m\}$ be the representatives of orbits of edges in $T$ under the action of $G$. The {\em transition matrix} $M_f$ is given by $[a_{ij}]_{1\leq i,j\leq m}$, where $a_{ij}$ is the number of times $f(e_j)$ crosses a translate of $e_i$. 
\begin{definition}\cite{bf,jpmutanguha}\label{defn_irred2}
	$f$ is {\bf $\GG$-irreducible} if $M_f$ is an irreducible matrix. 
\end{definition}
\begin{lemma}\cite[Lemma 3.6]{jpmutanguha}
	If $\phi$ is irreducible relative to $\GG$, then any $\GG$-relative weak representative of $\phi$ is $\GG$-irreducible.   
\end{lemma}
%\begin{proof}
%    Suppose $\phi$ is irreducible relative to $\GG$ but has a $\GG$-relative weak representative $f:T\to T$ that is $\GG$-reducible. Then $M_f$ is reducible. This implies that there exists $i,j\in\{1,\ldots,m\}$ such that the $(i,j)$-th element of $M_f^n$ for every $n>0$ is $0$. Then for no $n>0$, $f^n(e_j)$ crosses a translate of $e_i$, so there exists an $f$-invariant proper subforest of $T$, which does not contain any of the translates of $e_i$. This corresponds to an $f$-invariant free factor system $\{[G]\}\succ\GG'\succ\GG$.
%\end{proof}
\noindent
We also recall the equivalent definition of $\GG$-irreducibility in terms of a graph from \cite{bf}. Let $\Gamma$ be a finite oriented graph on $m$ vertices $\{x_1,\ldots,x_m$\}. There exist $a_{ij}$ many edges from $x_j$ to $x_i$. $M_f$ is irreducible if and only if for every $1\leq i,j\leq m$, there exists an oriented path from $x_i$ to $x_j$.

Now, for an irreducible $M_f$, there exists a unique eigenvalue $\lambda_f\geq 1$ and a corresponding eigenvector $u=(u_1,\ldots,u_m)$, and putting $l_T(e_i)=l_T(ge_i)=u_i$ for $1\leq i\leq m$ and any $g\in G$, we have the Perron-Frobeinus metric on $T$. We also refer to $\lambda_f$ as the stretch factor or growth rate of $f$. 
%\begin{lemma}\label{expanding}
%Let $f:T\to T$ be a train track map representing a fully irreducible non-surjective injective endomorphism $\phi:(G,\GG)\to(G,\GG)$. Suppose $T$ has at least two orbits of edges, with respect to the Perron-Frobeinus metric $d_T$ on $T$, then $\lambda_f>1$.
%\end{lemma}
%\noindent
%Clearly, $\lambda_f\geq 1$. The proof of $\lambda_f>1$ is the same as that of \cite[Lemma 1.11]{DahmaniLi} and therefore, we skip it. 
%\noindent
We will assume that any $\GG$-relative weak representative $f:T\to T$ is tight, i.e., for every $e\in E(T)$, $f$ is either injective in the interior of $e$ or $f(e)\in V(T)$. 

\noindent
For the rest of this subsection, let $\phi:(G,\GG)\to(G,\GG)$ be an injective endomorphism that is strictly-type preserving relative to $\GG$ and let $f:T\to T$ be a tight $\GG$-irreducible $\GG$-relative weak representative. The following results are from \cite[Section 1]{bf}. The results are very similar as there are only finitely many edge orbits for the action of $G$ on $T$. 
%\subsubsection{Collapse of a subforest}
%Let $Y\subset T$ be an $f$-invariant, $G$-invariant subforest with bounded components. If $e\in E(T)$ such that $f(e)\in V(T)$, $Y=\cup_{g\in G}\{ge\}$. 
%$f$ is a qi embedding, so there cannot be a collapsing subforest with unbounded components.
%Let $T_1=T/Y$ be obtained by collapsing each component of $Y$ to a point. This induces a map $T_1\to T_1$ and let $f_1:T_1\to T_1$ be obtained by tightening it. 
%\begin{lemma}\label{collapse}
%$f_1:T_1\to T_1$ is $\GG$-irreducible and the stretch factor $\lambda_{f_1}<\lambda_f$.  
%\end{lemma}
%\begin{proof}
%Let $\{e_{r+1},\ldots,e_m\}$ be the orbit representatives of edges in $Y$. Let $v=(v_1,\ldots,v_m)$ be the Perron-Frobenius eigenvector of $M_f$. Then $M_fv=\lambda_f v$. Let $M_1$ be the transition matrix of $f_1$. This is an $r\times r$ submatrix of $M_f$. The matrix $M_1$ is irreducible. Let $w=(v_1,\ldots,v_r)$. Then $M_1w=\lambda_{f_1}w$ and clearly $\lambda_{f_1}<\lambda_f$.  
%$M_fv=\lambda_f v$ implies for each $1\leq i\leq m$, $\sum_{1\leq j\leq m}a_{ij}v_j=\lambda_fv_i$ and $M_1w=\lambda_{f_1}w$ implies for each $1\leq i\leq r$, $\sum_{1\leq j\leq r}a_{ij}v_j=\lambda_1v_i<\lambda_fv_i$.
%\end{proof}
\subsubsection{Subdivision}
Let $w$ be a point in the interior of the edge $e_m$ in $T$ such that $f(w)\in V(T)$. Then we add $w$ to the vertex set of the tree, i.e., we have $T_1$ with $V(T_1)=V(T)\cup\bigcup_{g\in G}\{gw\}$ and $E(T_1)$ is obtained by replacing $ge_m$ by $ge'_m$ and $ge'_{m+1}$  such that $t(ge'_m)=gw=o(ge'_{m+1})$, for each $g\in G$. Let $f_1:T_1\to T_1$ be the new map.
\begin{lemma}\label{subdivide}
	$f_1:T_1\to T_1$ is a $\GG$-irreducible $\GG$-relative weak representative and $\lambda_{f_1}=\lambda_f$.    
\end{lemma}
\noindent
The proof is similar to that of \cite[Lemma 1.10]{bf}.   
%\begin{proof}
%    Suppose $f(e_i)$ pass through a translate of $f(e_m)$ for $i\in\{1,\ldots,m\}$. Then $f_1(e_i)$ passes through $e'_m$ and $e'_{m+1}$. Then for the new matrix $M_1:=M_{f_1}=[b_{ij}]_{1\leq i,j\leq m+1}$, $b_{im}=b_{im+1}=a_{im}$. And for $1\leq i\leq m-1$, $a_{mi}=b_{mi}+b_{m+1i}$. Take $w=(v_1,\ldots,v_m,v_m)$.
%\end{proof}
\subsubsection{Folding}
Let $e_{m-1},e_m$ be edges incident on the same vertex $v$ with $f(e_{m-1})=f(e_m)$. Then for each $g\in G$ we fold $ge_{m-1},ge_m$ to get a new tree $T_1$ with $E(T_1)=\{E(T)\setminus\{\cup_{g\in G}(ge_{m-1}\cup ge_m)\}\}\bigcup_{g\in G}\{ge'_{m-1}\}$. 
%We have $V(T_1)=\{V(T)\setminus\{\cup_{g\in G}(t(ge_{m-1})\cup t(ge_m))\}\}\bigcup_{g\in G}\{t(ge'_{m-1})\}$. 
Let $f_1:T_1\to T_1$ be the induced map with $f_1(ge'_i)=f(ge_i)$ for all $g\in G$ and $1\leq i\leq m-1$. If $f_1$ is not a $\GG$-relative weak representative, let $f_2:T_2\to T_2$ be the map obtained by collapsing anymaximal $f_1$-invariant subforest with bounded components, and tightening. Now if $f_1$ is a $\GG$-relative weak representative, we take $f_2=f_1$. 

Suppose $f(e_{m-1})\neq f(e_m)$, but have the same initial segment, then we apply subdivision on $e_{m-1},e_m$ and then the above folding on the new pair of edges. In both cases, we have the following result.
\begin{lemma}\label{fold}
	The map $f_2$ is $\GG$-irreducible and the stretch factor $\lambda_{f_2}\leq\lambda$.
\end{lemma}
\noindent
The proof is similar to that of \cite[Lemma 1.15]{bf}. It is easy to see that collapsing an maximal $f_1$-invariant subforest with bounded components or pretrivial edges reduces the stretch factor.  
%\begin{proof}
%   Suppose $f(e_{m-1})=f(e_m)$. We have $M_f=[a_{ij}]_{1\leq i,j\leq m}$ and the new matrix $M_1=[b_{ij}]_{1\leq i,j\leq m-1}$. Then for $b_{m-1i}=a_{m-1i}+a_{mi}$ for $1\leq i\leq m-1$ and for $1\leq i\leq m-2$, $b_{im-1}=a_{im-1}$. Since $a_{im}=a_{im-1}$, taking $w=(v_1,\ldots,v_{m-2},v_{m-1}+v_m)$, we get $M_fv=M_1w$. Therefore if $f_1$ is tight and has no $f_1$-invariant subforests or pretrivial edges, $\lambda_f=\lambda_1$. Otherwise, we tighten $f_1$ and apply Lemma \ref{collapse}.
%\end{proof}
\subsubsection{Removing valence-two vertices}
Let $v\in V(T)$ be a valence-two vertex. Without loss of generality, let $t(e_{m-1})=v=o(e_m)$. Note that for each $g\in G$, $gv$ is a valence-two vertex with $ge_{m-1}$ and $ge_m$ incident on $gv$. Let $H_t:T\to T$ be a homotopy with support in $\bigcup_{g\in G}\{ge_{m-1}\cup ge_m\}$ where $H_0$ is the identity map and $H_1$ is a map that collapses $ge_m$ to $t(ge_m)$ and stretches $ge_{m-1}$ over $ge_m$ for each $g\in G$. Let $f':T\to T$ be obtained by tightening $H_1\circ f$. Note that $gv\notin Im(f')$ for any $g\in G$. \\
Let $T_2$ be a tree with $V(T_2)=V(T)\setminus\{\cup_{g\in G}gv\}$ and $E(T_2)=\{E(T)\setminus\cup_{g\in G}(ge_{m-1}\cup ge_m)\}\bigcup_{g\in G}ge_{m-1}'$. For each $g\in G$, $ge'_{m-1}$ is an edge joining $o(ge_{m-1})$ and $t(ge_m)$, replacing $ge_{m-1}\cup ge_m$. Let $f_2:T_2\to T_2$ be the induced map and let $f_1:T_1\to T_1$ be obtained after collapsing the maximal $f_2$-invariant subforest with bounded components.
\begin{lemma}\label{valence-two}
	Suppose $f_1:T_1\to T_1$ is a $\GG$-irreducible $\GG$-relative weak representative obtained as above. Let $M_1$ be the transition matrix for $f_1:T_1\to T_1$ and let $\lambda_1$ be the stretch factor of $f_1$. If for the Perron-Frobenius eigenvector $v=(v_1,\ldots,v_m)$ of $M_f$, $v_{m-1}\leq v_m$ then $\lambda_1\leq\lambda_f$ and in particular, if $v_{m-1}<v_m$ then $\lambda_1<\lambda_f$.
\end{lemma}
\noindent
The proof is similar to the proof of \cite[Lemma 1.13]{bf}.
%\begin{proof}
%   Let $M_f=[a_{ij}]_{\leq i,j\leq m}$. Suppose $f_2=f_1$, i.e., $T_2$ has no $f_2$-invariant $G$-invariant subforests with bounded components and $T_2$ has no pretrivial edges. Then let $M=M_{f_2}=[b_{ij}]_{\leq i,j\leq m-1}$. Then $b_{im-1}=a_{im-1}+a_{im}$ and $b_{m-1i}=a_{m-1i}$. We can rearrange $e_{m-1},e_m$ so that $v_{m-1}\leq v_m$.

%   Now, for each $1\leq i\leq m-1$, $\sum_{j=1}^{m-1}b_{ij}v_j=\sum_{j=1}^{m-2}a_{ij}v_j+(a_{im-1}+a_{im})v_{m-1}$.

%   If $v_{m-1}=v_m$, we get $\sum_{j=1}^{m-1}b_{ij}v_j=\lambda_fv_i$ and so, $\lambda_1=\lambda_f$. Else, $\sum_{j=1}^{m-1}b_{ij}v_j<\lambda_fv_i$ and so, $\lambda_1<\lambda_f$.
%If $f_2\neq f_1$, clearly $\lambda_1<\lambda_f$.  
%\end{proof}
\subsubsection{Removing valence-one vertices}
Let $v\in V(T)$ be a valence-one vertex. Without loss of generality, let $e_m$ be incident on $v$. Then for each $g\in G$, $gv$ is a valence-one vertex and $ge_m$ is incident on it. Let $H_t:T\to T$ be a homotopy where $H$ is the identity map and $H_1$ is a map that collapses $ge_m$ to $gv$ for each $g\in G$. Let $f':T\to T$ be obtained by tightening $H_1\circ f$. Note that $gv\notin Im(f')$ for any $g\in G$. \\
Let $T_2$ be a tree with $V(T_2)=V(T)\setminus\{\cup_{g\in G}gv\}$ and $E(T_2)=E(T)\setminus\{\cup_{g\in G}ge_m\}$. Let $f_2:T_2\to T_2$ be the induced map and let $f_1:T_1\to T_1$ be obtained after collapsing the maximal $f_2$-invariant subforest with bounded components.
\begin{lemma}\label{valence-one}
	Suppose $f_1:T_1\to T_1$ is an irreducible $\GG$-relative weak representative obtained as above. Let $M_1$ be the transition matrix for $f_1:T_1\to T_1$ and let $\lambda_1$ be the stretch factor of $f_1$. Then $\lambda_1<\lambda_f$.
\end{lemma}
\noindent
Proof is similar to that given in \cite[Lemma 1.11]{bf}.
\noindent
We use these above operations to get a $\GG$-relative representative of $\phi$ and $T\in\TT_{\GG}$ with a train track structure, i.e., each vertex of $T$ has at least two gates. \\
Let $v\in V(T)$ with only one gate. Note that $v$ is free. For a non-free vertex with stabilizer $H_i$, if $e$ is an edge incident on it, for any $h\in H_i$, $e,he$ are in distinct gates as the arc stabilizers are trivial (See \cite[Lemma 6.6]{FranMartino}). \\
Also note that $v$ will be a finite valence vertex, in particular, the valence of $v$ will be bounded above by $m$. Otherwise, one can find a pair of edges $e,ge$ incident on $v$ and by the argument above, they cannot be in the same gate. 
Let $\{e_{i_1},\ldots,e_{i_s}\}$ be the set of edges incident on $v$. For $1\leq j\leq s$, let $f(e_{i_j})=x_1x_2\ldots x_kx^j_{k+1}\ldots x^j_{N_j}$ with $x_i,x^j_l\in E(T)$. Here, for distinct $j$, $x^j_{k+1}$ are also distinct. For $1\leq j\leq s$, let $w_j$ be a point in $e_{i_j}$ such that $f(w_j)=t(x_k)$. We apply subdivision on all the translates of $e_{i_j}$ in $T$ to get $f_1:T_1\to T_1$. Let the new vertex $w_{i_j}$ divide $e_{i_j}$ into $e'_{i_j},e''_{i_j}$.  
Then we apply folding on the turns $(e'_{i_j},e'_{i_l})$ for distinct pairs of $j,l$ and their translates. This results in $v$ being a valence-one vertex. Then applying Lemma \ref{valence-one}, we get $f_2:T_2\to T_2$. \begin{lemma}\label{ttstructure}
	There exists $T\in\TT_{\GG}$ and a $\GG$-irreducible $\GG$-relative weak representative $f:T\to T$ of $\phi$ such that the gate stricture induced by $f$ is a train track structure.
\end{lemma}
\noindent
Further, we have the following.
%Since the valence of $v$, and hence of $gv$, $g\in G$, is finite, apply the above steps to all the edges emanating from $v$ and $gv$. Then applying Lemmas \ref{collapse}, \ref{subdivide} \ref{fold} and \ref{valence-one}, we get the following.
%\begin{lemma}\label{gates}
%Let $f_3:T_3\to T_3$ be a $\GG$-relative weak representative of $\phi$ with a train track structure. Then $f_3$ is $\GG$-irreducible and the stretch factor $\lambda_3\leq\lambda_f$.   
%\end{lemma}

\begin{theorem}\cite[Proposition 3.7]{jpmutanguha}\label{minimal-stretch}
	Let $\phi:(G,\GG)\to(G,\GG)$ be an injective endomorphism that is irreducible relative to $\GG$ and strictly type-preserving relative to $\GG$. Then there exists a $\GG$-relative representative $f:T\to T$ of $\phi$ with minimal stretch factor.
\end{theorem}
%\begin{proof}
%    Let $T\in\TT_{\GG}$. The quotient of $T$ under the action of $G$ has at most $N:=p+3k-3$ orbits of edges. Then the irreducible transition matrix of a $\GG$-relative representative $g:T\to T$ is an $n\times n$-matrix, where $n\geq N$. Let $B$ be an irreducible non-negative $n\times n$ integer matrix with the Perron Frobenius eigenvalue $\lambda\leq\lambda_g$. Let $v$ be the Perron Frobenius eigenvector and without loss of generality, assume that the smallest entry of $v$ is $1$. Then the minimum row sum of $B^k$, $k\geq1$ is bounded above by $\lambda^k$. Also, the largest entry of $B$ is bounded above by the minimum row sum of $B^2$. Thus, every element of $B$ is bounded by $\lambda_g^2$. So, number of candidates for $B$ is finite and the number of possible stretch factors is finite. This proves the theorem.
%\end{proof}
%\noindent
%With slight modification, the proof of \cite[Proposition 3.7]{jpmutanguha} works in this case as well. 
\begin{theorem}\cite[Theorem 1.7]{bf}\label{traintrack}
	Let $\phi:(G,\GG)\to(G,\GG)$ be an injective endomorphism that is irreducible relative to $\GG$ and strictly type-preserving relative to $\GG$. Then there exists a $\GG$-relative representative of $\phi$ that is a train track map.
\end{theorem}
\noindent
{\it Sketch of the proof:}
Let $f:T\to T$ be a a $\GG$-relative representative of $\phi$ with the minimal stretch factor $\lambda$. By minimality of the stretch factor, $T$ has no valence-one vertices and by Lemma \ref{valence-two}, we may assume that $T$ has no valence-two vertices. By Lemma \ref{ttstructure}, we may also assume that the gate structure on $T$ induced by $f$ is a train track structure. 

Suppose there exists a legal turn $(e,e')$ at $v\in V(T)$ that is mapped to an illegal turn at $f(v)\in V(T)$. Let $f^2(e)=x_1x_2\ldots x_kx_{k+1}\ldots x_N$ and $f^2(e')=x_1x_2\ldots x_kx'_{k+1}\ldots x'_M$, with $x_i,x'_j\in E(T)$ for $1\leq i\leq N, k+1\leq j\leq M$. Let $w,w'$ be the points in $e,e'$ respectively such that $f^2(w)=t(x_k)=f^2(w')$. We apply subdivision and get $w,w'$ and their translates as new vertices. Note that $w,w'$ are valence-two. We fold $(e,e')$ by Lemma \ref{valence-two}. We repeat this process for every $G$-translate of $(e,e')$, and also to every legal turn in $T$ that is mapped to an illegal turn.  Let $f_1:T_1\to T_1$ denote the new map. Subdivision does not change the stretch factor and by the minimality of $\lambda_f$, removing the valence-two vertex also does not change the stretch factor. So every legal turn in $T_1$ is mapped to a legal turn under $f_1$ and the stretch factor of $f_1$ is $\lambda$. Note that $T_1$ cannot have valence-one vertices. Thus, condition $(1)$ is satisfied. 

Now by the above arguments, suppose $f$ satisfies condition $(1)$. When $\lambda_f=1$, $f(e)$ is an edge for every edge $e$ in $T$ and so $f(e)$ is legal. So let $\lambda_f>1$. Let $e\in E(T)$ (and its translates) such that $f(e)$ is not a legal path. Let $w$ be the point in the interior of $e$ which fails to be locally injective under $f^2$. We apply subdivision and get $f_1:T_1\to T_1$ with $gw\in V(T_1)$ for all $g\in G$ and with $ge$ divided into two edges for each $g\in G$. Here, we have $f_1:T_1\to T_1$ where a legal turn at $w$ is mapped to an illegal turn at $f_1(w)$. We apply the above steps to get $f_2:T_2\to T_2$ which satisfies condition $(2)$ and therefore, is a train track map.
\qed

\vspace{2mm}
\noindent
Let $k>1$. Since $\phi$ is fully irreducible relative to $\GG$, by Theorem \ref{traintrack}, $\phi^k$ has a train track representative. Clearly, $f^k$ is a $\GG$-relative weak representative of $\phi^k$. Suppose there exists a legal turn $\tau_0=(e,e')$ that is mapped to an illegal turn $\tau_1$ by $f^k$ and for every $l<k$, it is mapped to a legal turn under $f^l$. In particular, $f^{k-1}$ maps $\tau_0$ to a legal turn $\tau$ and $f$ maps $\tau$ to $\tau_1$. This cannot happen as $f$ is a train track map. Similarly, if $f^k(e)$ is not a legal path, then there exists an illegal turn $\tau_1$ in $f^k(e)$. Similarly, $f^k(e)$ is a legal path for any $e\in E(T)$. Thus, we have the following. \begin{lemma}\label{power-traintrack}
	Let $\phi:(G,\GG)\to(G,\GG)$ be an injective endomorphism that is fully irreducible relative to $\GG$ and strictly type-preserving relative to $\GG$. If $f:T\to T$ is a train track map representing $\phi$, then for any $k>1$, $f^k:T\to T$ is a train track map representing $\phi^k$.
\end{lemma}
\noindent
Now we recall the following definition. 
\begin{definition}\cite{DahmaniLi} 
	Let $\PP$ be a $\phi$-invariant collection of conjugacy classes of subgroups of $G$. We say $\PP$ is \textbf{hyperbolically coning-off $\GG$} if\\ (1) $\PP$ is a finite collection of conjugacy classes of subgroups of $G$,\\
	(2) $\PP$ is malnormal,\\
	(3) For every $[H_i]\in\GG$, there exists $[P_j]\in\PP, g\in G$ such that $H_i\leq gP_jg^{-1}$, and\\
	(4) for all $i$, the action of $P_i$ on $T_{P_i}$ is cofinite.
\end{definition}
\noindent
Let $\PP=\{[P_1],\ldots,[P_q]\}$ for the rest of this paper. For each $1\leq i\leq q$, let $T_i$ denote the minimal subtree of $P_i$ in $T$. Let $\widehat{T}=\EE(T,\{T_i\}_{i=1}^q)$ be the tree obtained by coning-off $T_i's$ and their translates in $T$. 
\begin{lemma}\cite[Proposition 1.12]{DahmaniLi}\label{relhypG}
	$G$ is strongly hyperbolic relative to $\PP$. 
\end{lemma}
\noindent
Given a train track map $f:T\to T$, this induces a map $\hat{f}:\widehat{T}\to\widehat{T}$ defined as follows: for $[P]\in\PP$, there exists $[P']\in\PP,g\in G$ such that $\phi(P)\leq gP'g^{-1}$. Let $v_P$ and $v_{gP'}$ denote the corresponding cone points. Then $\hat{f}(\nu_{P})=\nu_{gP'}$. For the rest of the paper, we assume that every turn at a cone point is legal. Then, if $f:T\to T$ is a train track map, so is $\hat{f}:\widehat{T}\to\widehat{T}$.

%\begin{proof}
%    Let $e\in E(\widehat{T})$. If $e\in E(T)$, we know that $f(e)$ is a legal path. So let $e$ be an edge joining a cone point $\nu_{P}$ and some element $h\in P$, where $[P]\in\PP$. Then, by the definition of the map $\hat{f}$, $\hat{f}(e)$ is an edge joining a cone point to an element in $G$. Thus, $\hat{f}$ maps any edge to a legal path.

%    Now, for any $v\in V(T)$, a legal turn at $v$ is mapped to a legal turn at $\hat{f}(v)$. So, it is enough to check for legal turns at any cone-point. \textcolor{red}{to be completed}
%\end{proof}
\noindent
Given $(G,\GG), \phi:(G,\GG)\to(G,\GG)$ fully irreducible relative to $\GG$, $T\in\TT_{\GG}$, $\PP$ as above. 
\begin{lemma}[Bounded cancellation lemma]\label{bcl}\cite[Lemma 2.9]{DahmaniLi} Let $f:T\to T$ be a $\GG$-relative weak representative of $\phi$.
	There exists $C_{\ref{bcl}}>0$ depending on $f,\PP$ such that for any path $\rho=\alpha\cdot\beta$ (without cancellation) in $T$, we have $f(\rho)=\alpha_1\cdot\beta_1$, where $f(\alpha)=\alpha_1\cdot u$ and $f(\beta)=\Bar{u}\cdot\beta_1$ such that $u$ is a path in $T$ with $l_T(u)\leq C_{\ref{bcl}}$. Moreover, we have $$l_T([f([\rho])])\geq l_T([f([\alpha])])+l_T([f([\beta])])-2C_{\ref{bcl}}.$$  
\end{lemma}
\begin{remark}\begin{enumerate}
		\item By \cite[Lemma 2.9]{DahmaniLi}, we also have $$l_{\T}([\f([\hrho])])\geq l_{\T}([\f([\hat{\alpha}])])+l_{\T}([\f([\hat{\beta}])])-2C_{\ref{bcl}}.$$
		\item By induction, if $\rho=\alpha_1\cdot\alpha_2\cdots\alpha_m$ (without cancellation) in $T$, then $$l_T([f([\rho])])\geq \sum_{i=1}^ml_T([f([\alpha_i])])-2(m-1)C_{\ref{bcl}}.$$
	\end{enumerate}
\end{remark}
\noindent
Now, we have the following results about indivisible Nielsen paths.
\begin{lemma} 
	An indivisible Nielsen path $\rho$ has a unique illegal turn.
\end{lemma}
\begin{proof}
	Let $\rho$ be a geodesic. Let $N>0, g\in G$ such that $[f^N(\rho)]=g\rho$. Suppose   
	$\rho=\rho_1\cdot\ldots\cdot\rho_k$, where $k\geq 3$ and each $\rho_i$ is a maximal legal segment, i.e., $\rho$ has $k-1$ illegal turns. By bounded cancellation lemma, $[f^N(\rho)]=g\rho_1\cdot\ldots\cdot g\rho_k$, where $[f^N(\rho_1)]=g\rho_1u_1$, $[f^N(\rho_i)]=\bar{u}_{i-1}g\rho_i u_i$, for $2\leq i\leq k-1$ and $[f^N(\rho_k)]=\bar{u}_{k-1}g\rho_k$. 
	
	Put $\rho^0_i=\rho_i$, $2\leq i\leq k-1$. We have $f^N(\rho_i)=\bar{u}_{i-1}g\rho_i u_i$. Let $\rho_i=\alpha_1\cdot\rho^1_i\cdot\beta_1$ such that $f^N(\rho^1_i)=g\rho_i=g\rho^0_i$.
	Inductively, for $j\geq1$, let $\rho_i^j=\alpha_{j+1}\cdot\rho^{j+1}_i\cdot\beta_{j+1}$, such that $f^N(\rho^{j+1}_i)=g\rho^j_i$. Let $x_i=\bigcap_{j\geq0}\rho^j_i$. Then $f^N(x_i)=gx_i$, for $2\leq i\leq k-1$. This contradicts the indivisibility of $\rho$.
\end{proof}
\begin{lemma}\label{bdd-nielsen}
	The number of indivisible Nielsen paths is uniformly bounded, up to their $G$-orbits.   
\end{lemma}
\begin{proof}
	Let a geodesic $\rho$ be an indivisible Nielsen path of period $n$. Then $[f^n(\rho)]=g\rho$. Let $\rho=\alpha\cdot\beta$, where $\alpha,\beta$ are maximal legal segments. Then by bounded cancellation lemma, we have $[f^n(\alpha)]=g\alpha\cdot u$ and $[f^n(\beta)]=\bar{u}g\beta$, where $l_T(u)$ is uniformly bounded. Then $$l_T([f^n(\rho)])=\lambda^n l_T(\rho)-2l_T(u).$$
	This implies, $$(\lambda^n-1)l_T(\rho)=2l_T(u)\leq 2C_{\ref{bcl}}.$$ Thus, $l_T(\rho)\leq \frac{2C_{\ref{bcl}}}{\lambda^n-1}$. So, for any indivisible Nielsen path $\rho$, we have $l_T(\rho)\leq 2C_{\ref{bcl}}$, and as the action of $G$ on $T$ is cofinite, this proves the lemma. \end{proof}
\begin{remark}\label{nielsen}\begin{enumerate}
		\item Let $\alpha,\beta$ be pre-Nielsen paths and let $\rho=\alpha\cdot\beta$. Then there exists $N>0,g,h\in G$ such that $[f^N(\alpha)]=g\alpha$ and $[f^N(\beta)]=h\beta$. Since the terminal vertex of $\alpha$ and the initial vertex of $\beta$ are the same, $g=h$ and thus, $\rho$ is also pre-Nielsen. 
		\item Lemma \ref{bdd-nielsen} and (1) imply that there exists $M_{\ref{nielsen}}\geq1$ such that it is impossible to concatenate more than $M_{\ref{nielsen}}$ indivisible pre-Nielsen paths. 
\end{enumerate}\end{remark}
\noindent
For a path $\rho$ in $T$, we denote the number of illegal turns by $i(\rho)$.
\begin{lemma}\cite[Lemma 2.13]{DahmaniLi}\label{trichotomy}
	Let $f:T\to T$ be a train track map. Given $C>0$, let $M>0$ such that for any path $\rho$ in $T$, if the length of the longest legal segment in $f^M(\rho)$ is less than $C$, then one of the following happens.\\
	(i) $i(f^M(\rho))< i(\rho)$.\\
	(ii) $\rho=\gamma_1\cdot\alpha_1\cdots\alpha_s\cdot\gamma_2$, where $\gamma_i,\alpha_j$, $i=1,2,1\leq j\leq s$ are pre-Nielsen paths with at most one illegal turn each and length of $\gamma_1,\gamma_2$ is at most $2C$.
\end{lemma}
%\begin{proof}
%Suppose (i) does not hold. Then $i(f^n(\rho))=i(\rho)$ for all $n>0$. Note that $l_T(f^n(\rho))\leq C\cdot (i(\rho)+1)$.\\
%Case 1: Angles in $[f^n(\rho)]$ are bounded for $n>0$.\\
%Then, there exists $N>0, n\geq0, g\in G$ such that $f^{N+n}(\rho)=gf^N(\rho)$. This implies that $\rho$ is a pre-Nielsen path. Let $\rho=\rho_1\cdot\rho_2\cdot\ldots\cdot\rho_k$, where each $\rho_i$ is a maximal legal segment. Then, for $2\leq i\leq k-1$, there exists $\tau_i\subset f^N(\rho_i)$ such that $g\tau_i\subset f^n(\tau_i)$ and so, there exists $x_i$ in the interior of $\rho_i$ satisfying $f^{N+n}(x_i)=gf^N(x_i)$. We can express $\rho=\gamma_1\cdot\alpha_1\cdot\ldots\cdot\alpha_s\cdot\gamma_2$, where each $\gamma_i$ and $\alpha_i$ is an indivisible pre-Nielsen path. 

%\noindent
%Case 2: $\{f^n(\rho)\}$ has unbounded angles. By Lemma \ref{unbdd-angle}, there exists a vertex, say $v$ in $\rho$ such that as $n\to\infty$, angle at $\f^n(v)$ tends to infinity. Let $v_1,v_2$ be the end points of $\rho$ and let $\rho=\rho_1\cdot\rho_2$, where $\rho_i$ is the subpath of $\rho$ joining $v_i,v$, $i=1,2$. Then angles in $[f^n(\rho_i)]$ is bounded. As in Case 1, we have that $\rho_i$ can be expressed as a concatenation of indivisible pre-Nielsen paths. Thus, $\rho$ is a concatenation of indivisible Nielsen paths.
%\end{proof}
\subsection{Polynomially growing subgroups}
Let $G=H_1\ast\cdots\ast H_p\ast F_k$, where $\GG=\{[H_1],\ldots,[H_p]\}$ and let $\phi:(G,\GG)\to(G,\GG)$ be an injective endomorphism that is fully irreducible relative to $\GG$. Let $f:T\to T$ be a train track map representing $\phi$, for some $T\in\TT_{\GG}$. For an element $g\in G$, the
length $||g||$ of its conjugacy class is the length of a cyclically reduced element conjugate to $g$. 
\begin{definition}
	An element $g\in G$ is \textbf{polynomially growing for $\phi$} if $||\phi^n(g)||$ is bounded by a polynomial in $n$ for all $n\geq0$.
	
	A subgroup $H\leq G$ is \textbf{polynomially growing for $\phi$} if every element of $H$ is polynomially growing for $\phi$.
\end{definition} 
\noindent
Let $T\in\TT_{\GG}$. For each $n\geq1$, let $$d_n(x,y)=\frac{d_T(f^n(x),f^n(y))}{\lambda_f^n}$$ for every $x,y\in T$. Let $d_{\infty}(x,y)=\lim_{n\to\infty}d_n(x,y)$ and let $T_{\infty}=T/\sim$ where $x\sim x'$ if $d_{\infty}(x,x')=0$. 
It is easy to see that $(T_{\infty},d_{\infty})$ is an $\RR$-tree. Also, $f$ induces a map $F:T_{\infty}\to T_{\infty}$ satisfying $d_{\infty}(F(x),F(y))=\lambda_f\cdot d_{\infty}(x,y)$ for all $x,y\in T_{\infty}$.
We recall the following definition from \cite{horbez}.
\begin{definition}
	A $(G,\GG)$-tree is {\bf small} if arc stabilizers are either trivial or cyclic.
	It is {\bf very small} if it is small, the arc stabilizers are closed under taking roots, and not conjugate into any of the elements in $\GG$ and tripod stabilizers are trivial.
\end{definition}
%\begin{proposition}\cite[Proposition 3.1]{horbez}\label{horbez}
% The space of projective very small, minimal, nontrivial $\RR$-trees such that for each $[H_i]\in\GG$, $H_i$ fixes a point of in  is closed.   
%\end{proposition}
\begin{lemma}\cite[Lemma 2.7]{gaboriau}\label{non-trivial}
	Action of $G$ on $T_{\infty}$ is non-trivial.  
\end{lemma}
\begin{proof}
	Let $\{e_1,\ldots,e_m\}$ be the edge orbit representatives of action of $G$ on $T$ and let $f_i$ denote the image of $e_i$ in $T_{\infty}$, $1\leq i\leq m$. Now suppose the action is trivial, i.e., there exists $x_0\in T_{\infty}$ fixed by $G$. Without loss of generality, assume $x_0\in f_1$. Also, let the diameter of the smallest finite subtree of $T_{\infty}$ containing $f_1,\ldots,f_m$ be $D$. Note that for any $y\in T_{\infty}$, $y\in g\cdot f_i$, for some $g\in G, i\in\{1,\ldots,m\}$, so $y$ lies in a finite neighbourhood of $x_0$, i.e., $T_{\infty}$ is bounded. This is a contradiction since the length of $f_i$'s increase exponentially under iterations of $F:T_{\infty}\to T_{\infty}$.
\end{proof}
\noindent
\begin{lemma}\cite[Lemma 2.8]{gaboriau}\label{arc-stab-trivial} $T_{\infty}$ has trivial arc stabilizers.
\end{lemma} 
\noindent
We give a slightly different proof for this. We first recall the following results. For an $\RR$-tree $T$ and $g\in Isom(T)$, let $C_g=\{x\in T\mid d_T(x,gx)=l_T(g)\}$ be the characteristic set of $g$.
\begin{lemma}\cite[Lemma 5.2]{CM}\label{5.2}
	Let $p_1,p_2$ be points on an $\RR$-tree $T$ and $g_1,g_2\in Isom(T)$ such that $d_T(p_i,g_i\cdot p_i)<\epsilon$, for $i=1,2$. If $d_T(p_1,p_2)>2\epsilon$ then $C_{g_1}\cap C_{g_2}$ contains a geodesic of length at least $d_T(p_1,p_2)-2\epsilon$.
\end{lemma}
\begin{lemma}\cite[Lemma 1.10]{CM}\label{1.10}
	Suppose $g,h\in Isom(T)$ and $C_g\cap C_h$ contains an arc $J$ with length $l(J)\geq l_T(g)+l_T(h)$. If $g,h$ are both hyperbolic then also assume that the orientations induced by $g,h$ on their axes agree on the overlap. Then $ghg^{-1}h^{-1}$ fixes a subarc of $J$ of length $l(J)-l_T(g)-l_T(h)$.
\end{lemma}
\noindent
\textit{Proof of Lemma \ref{arc-stab-trivial}:}
Suppose there exists an arc $c$ in $T_{\infty}$ and $g'\in G\setminus\{1\}$ such that $g'\cdot c=c$. Choose $k>0$ large enough so that there exists subarcs $c_0,c_1$ of $F^k(c)$ and $h'\in G\setminus\{1\}$ such that $h'\cdot c_0=c_1$. Put $g=\phi^k(g')$ and $h=h'^{-1}gh'$. Then $g\cdot c_0=c_0$ an $h\cdot c_0=c_0$. Let $x,y$ be the endpoints of $c_0$ in $T_{\infty}$. We make the following observations. \\
Let $d_{\infty}(x,y)=L$. Then $\lim_{n\to\infty}d_n(x,y)=L$ and for each $\epsilon>0$, there exists $N_{\epsilon}>0$ such that for each $n\geq N_{\epsilon}$, we have $L-\epsilon< d_n(x,y)<L+\epsilon$.\\
Since $d_{\infty}(x,g\cdot x)=0=d_{\infty}(y,h\cdot y)$, we have that for each $\epsilon>0$, there exists $N'_{\epsilon}>0$ such that $d_n(x,g\cdot x)<\epsilon$ and $d_n(y,h\cdot y)<\epsilon$ for each $n\geq N'_{\epsilon}$.\\
Also, $l_{T_{\infty}}(g)\leq d_{\infty}(x,gx)=0$. Similarly, $l_{T_{\infty}}(h)\leq d_{\infty}(y,hy)=0$. Then for a given $\epsilon>0$, there exists $N''_{\epsilon}>0$ such that $l_{T_n}(g)\leq d_n(x,gx)<\epsilon$ for all $n\geq N''_{\epsilon}$. Similarly, for all $n\geq N''_{\epsilon}$, $l_{T_n}(h)<\epsilon$.

Choose $\epsilon>0$ such that $\epsilon<\frac{L}{6}$. Then for $N=\max\{N_{\epsilon},N'_{\epsilon},N''_{\epsilon}\}$, $d_N(x,y)>5\epsilon$ and by Lemma \ref{5.2}, there exists an arc $\alpha$ in $C_g\cap C_h$ in $T_N$ such that $l_N(\alpha)>3\epsilon$. We have $l_N(\alpha)>3\epsilon>l_{T_N}(g)+l_{T_N}(h)$. So by Lemma \ref{1.10}, there exists a subarc $\beta$ of $\alpha$ with length $l_N(\beta)\geq\epsilon$ that is fixed by $ghg^{-1}h^{-1}$, without loss of generality. This is a contradiction since $G$ acts on $T_N$ with trivial arc stabilizers.
\qed

\vspace{2mm}

\noindent
Recall that for $G=H_1\ast\cdots\ast H_k\ast\FF_p$ with $\GG=\{[H_1],\ldots,[H_k]\}$, the Kurosh rank of $(G,\GG)$ is $k+p$. By \cite[Theorem I]{MCML}, the space of projective, very small, minimal, non-trivial $G$-trees is compact and thus, $T_{\infty}$ is a very small tree. 
\begin{lemma}\cite[Corollary 5.5]{horbez}\label{kurosh-rk}
	(1) A very small $(G,\GG)$-tree has at most $2(k+p)-2$ orbits
	of branch or inversion points.\\
	(2) For a very small $(G,\GG)$-tree $T$, Kurosh rank of the stabilizer of any $x\in T$ is at most the Kurosh rank of $(G,\GG)$.
\end{lemma}
\noindent
With these results, the following result is true in our case as well. 
\begin{proposition}\cite[Proposition 1.13]{DahmaniLi}\label{polynomial}
	For a polynomially growing subgroup $P_0$ of $G$ for $\phi$ on $T$,
	(resp. for a polynomially growing element $h$ for $\phi$ in $T$), \begin{enumerate}
		\item there exists a unique maximal polynomially growing subgroup $P$ of $G$ for $\phi$ on $T$, such that $P_0\leq P$ (resp. such that $h\in P$).
		\item The collection of maximal polynomially growing subgroup $\PP$ of $G$ for $\phi$ on $T$ is malnormal.
		\item Any maximal polynomially growing subgroup of $G$ for $\phi$ on $T$ acts cofinitely on its minimal subtree in $T$ and has Kurosh rank smaller than the rank of $G$ (strictly if it is not $G$).
		\item If $P$ and $P'$ are two maximal polynomially growing subgroup of $G$ for $\phi$ on $T$, such that there exists $k\geq 1$ for which $\phi^k(P)\leq P$ and $\phi^k(P')\leq P'$ then $P=P'$.
		\item There are only finitely many conjugacy classes of maximal polynomially growing subgroups for $T$ in $G$. \end{enumerate}
\end{proposition}
%The proof is by induction on Kurosh rank of $G$. Without loss of generality, assume $G$ is not polynomially growing. If $P_0$ (resp $h$) is polynomially growing for $\phi$, then it stabilizes a point $x\in T_f$. Let $G_0<G$ be the stabilizer of $x$ in $T_{\infty}$.  Then $P_0\leq G_0$ (resp $h\in G_0$). By Lemma \ref{kurosh-rk}, the Kurosh rank of $G_0$ is strictly less than that of $G$. Then (1) follows easily. Let $P_1,P_2\in\PP$. Let $G_1,G_2$ be point stabilizers in $T_{\infty}$ such that $P_i\leq G_i$, for $i=1,2$. If $G_1\neq G_2$, then $P_1, P_2$ are malnormal as arc stabilizers are trivial. If $G_1=G_2$, the result follows by induction. This proves (2). 

%For (3), $G_0$ acts on $T$ with trivial arc stabilizers since $G_0<G$. Now, Kurosh rank of $G_0$ is finite. So, the action of $G_0$ on its minimal subtree is cofinite. By induction hypothesis, the action of a maximal polynomially growing subgroup $P\leq G_0$ on its minimal subtree in $T_{G_0}$ is cofinite. But the minimal subtree of $P$ in $T_{G_0}$ is the minimal subtree of $P$ in $T$. Thus we are done. The result of Kurosh rank of $P$ follows by induction.

%(4) holds since $d_{\infty}(F^k(x),F^k(y))=\lambda^k d_{\infty}(x,y)$ for all $x,y\in T_{\infty}$. By induction hypothesis, there are only finitely many conjugacy classes of maximal polynomially growing subgroups for $T$ in $G_0$ and there are only finitely many conjugacy classes of $G_0$ in $G$, thus (5) follows.
\noindent
Note that $\PP$ hyperbolically cones-off $\GG$.
\begin{lemma}\label{transversality}\cite[Proposition 1.15]{DahmaniLi}
	Suppose $\PP$ is a collection of conjugacy classes of maximal polynomially growing subgroups for $\phi$ on $T$. There exists $C_{\ref{transversality}}>1$ such that if $\rho$ is a path in $T$, then $l_{\T}(\hrho)\geq\frac{1}{C_{\ref{transversality}}}l_T(\rho)$.
\end{lemma}
\noindent
$C_{\ref{transversality}}$ is called the \textit{tranversality constant of legal paths with respect to $\PP$}.
\begin{definition}
	\textbf{Critical constant} for $f$, $C(f):=\frac{2C_{\ref{bcl}}}{\frac{\lambda_f}{C_{\ref{transversality}}}-1}$. 
\end{definition}
\begin{remark}
	Up to taking a power of $f$, we assume $\lambda_f>\max\{C_{\ref{transversality}},2C_{\ref{bcl}}+1\}$.
\end{remark}
\noindent
We recall the following corollary of the bounded cancellation lemma. 
\begin{lemma}\label{bcl-appl}\cite[Lemma 2.10]{DahmaniLi}\cite[Lemma 5.2]{brinkmann}
	Let $\beta$ be a legal segment in $T$ with $l_{\T}([\beta])\geq C(f)$. Let $\alpha,\gamma$ be paths in $T$ such that the concatenation $\alpha\cdot\beta\cdot\gamma$ is locally injective. Then there exists a constant $\nu_{\ref{bcl-appl}}>0$ such that the length of the segment in $\f^i([\alpha\cdot\beta\cdot\gamma])$ corresponding to $\beta$ is at least $\nu_{\ref{bcl-appl}}(\frac{\lambda_f}{C_{\ref{transversality}}})^il_{\T}(\beta)$ for all $i\geq 1$. 
\end{lemma}
\noindent
The following result is a trivial generalization of \cite[Lemma 5.3]{brinkmann}. 
\begin{lemma}\label{illegal_turns}
	Let $\rho$ be a path such for some $n\geq1$, any legal segment of $f^n(\rho)$ has length less than $C(f)$ and $i(f^n(\rho))\geq M_{\ref{nielsen}}$. Then, for all $m\geq1$, $$i(\rho)>\big(\frac{M_{\ref{nielsen}}+1}{M_{\ref{nielsen}}+2}\big)^mi(f^{nm}(\rho)).$$ 
\end{lemma}
\begin{proof}
	Let $f^n(\rho)=\rho_1\cdot\rho_2\cdots\rho_m\cdot\tau$ be a decomposition such that $i(\rho_i)=M_{\ref{nielsen}}$, for $1\leq i\leq m$ and $i(\tau)<M_{\ref{nielsen}}$. Let $\rho'_1,\ldots,\rho'_m,\tau'$ be subpaths of $\rho$ such that $f^n(\rho'_i)=\rho_i$, $1\leq i\leq m$, and $f^n(\tau')=\tau$. Note that for every $1\leq i\leq m$, the length of the maximal legal subsegment in $\rho'_i$ is less than $C(f)$. 
	%Otherwise, as the length of the legal segments grow exponentially, it will contradict the hypothesis that $f^n(\rho)$ has no legal segments of length greater than $C_f$. 
	Suppose $i(\rho'_i)=M_{\ref{nielsen}}$. By Lemma \ref{trichotomy}, $\rho'_i=\gamma_1\cdot\alpha_1\cdots\alpha_{M_{\ref{nielsen}}-2}\cdot\gamma_2$, where each $\gamma_j,\alpha_k$ is a pre-Nielsen path with one illegal turn and $\gamma_1,\gamma_2$ also have an illegal turn each. This contradicts Remark \ref{nielsen}(2). Therefore by Lemma \ref{trichotomy}, for $1\leq i\leq m$, we have $i(\rho'_i)\geq M_{\ref{nielsen}}+1$ and,
	
	\noindent
	$i(\rho)\geq \sum_{1=i}^{m}i(\rho'_i)+i(\tau')\geq (M_{\ref{nielsen}}+1)m+i(\tau')>M_{\ref{nielsen}}m+i(\tau)>\frac{M_{\ref{nielsen}}+2}{M_{\ref{nielsen}}+1}i(f^n(\rho))$.\\
	By induction, for all $k\geq1$, $$i(\rho)>\big(\frac{M_{\ref{nielsen}}+1}{M_{\ref{nielsen}}+2}\big)^ki(f^{nk}(\rho)).$$
\end{proof}
%\begin{remark}\label{illegal turns-coned}
%   Replacing $\rho$ by $\hrho$, for any $n\geq1$, if every legal segment of $\f^n(\hrho)$ is of length less than $C(f)$ and $i(\f^n(\hrho))\geq M_{\ref{nielsen}}$, we get $i(\hrho)>\big(\frac{M_{\ref{nielsen}}+1}{M_{\ref{nielsen}}+2}\big)^mi(\f^{nm}(\hrho))$.
%\end{remark}
\begin{lemma}\cite[Lemma 5.4]{brinkmann}\label{length-illegal}
	Given $C>0$, there exists $K_{\ref{length-illegal}}\geq 1$ such that for a path $\rho$ with each maximal legal segment of length less than $C(f)$ and $i(\rho)>0$, we have $$
	K_{\ref{length-illegal}}^{-1}\cdot i(\rho)\leq l(\rho)\leq K_{\ref{length-illegal}}\cdot i(\rho).$$
\end{lemma}

	%\section{Combination theorem}\label{section 4}
	%\input combination.tex
	
	%\input section4-240620.tex
	
	\section{The combination theorem}\label{section 5}
	\noindent
Given, $G=H_1\ast\cdots\ast H_p\ast\FF_k$, where $\GG=\{[H_1],\ldots,[H_p]\}$. Let $\phi:(G,\GG)\to(G,\GG)$ be an injective, non-surjective endomorphism that is fully irreducible relative to $\GG$ and strictly-type preserving relative to $\GG$. Let $\PP=\{[P_1],\ldots,[P_q]\}$ be $\phi$-invariant collection of conjugacy classes of subgroups of $G$ that hyperbolically-cones off $\GG$.\\
The group $G{\ast}_{\phi}$ is a finite graph of groups where the underlying graph is a loop with the vertex group $G$ and the edge group isomorphic to $G$ with identity and $\phi$ as the attaching maps. Let $T\in\TT_{\GG}$ and $f:T\to T$ be the train track map of $\phi$. Let $p:X\to\mathcal{T}$ be the corresponding tree of relatively hyperbolic spaces and $p:\mathcal{TC}(X)\to\mathcal{T}$ be the induced tree of coned-off spaces with copies of $\widehat{T}=\EE(T,\{T_{P}\}_{[P]\in\PP})$ as vertex and edge spaces. We denote the map induced on coned-off trees by $\hat{f}:\T\to\T$. Note that qi-embedded and qi-preserving electrocution conditions are satisfied for $G{\ast}_{\phi}$. Let $p:X\to\mathcal{T}$ be the corresponding tree of relatively hyperbolic spaces that is strictly-type preserving and satisfies the qi-embedded condition and the qi-preserving electrocution condition.\begin{definition}\cite{mit0,MjReeves}. 
	\begin{enumerate}
		\item $H:[-m,m]\times I\to X$ is a {\em hallway} of length $2m$ if:\\
		$(i)$ $H^{-1}(\bigcup_{v\in V(\TT)} X_v)=\{-m,\ldots,m\}\times I$,\\
		$(ii)$ $H$ maps $\{i\}\times I$ to a geodesic in some $X_v$.\\
		$(iii)$ $H$ is transverse, relative condition (i), to $\bigcup_e X_e$. 
		\item For $b>0$, a hallway is {\em $b$-thin} if $d(H(i,t),(i+1,t))\leq b$ for all $-m\leq i\leq m, t\in I$. It is {\em essential} if $p(H([-m,m]\times I)$ does not backtrack.
		\item For $\lambda>1$, a hallway of length $2m$ is {\em $\lambda$-hyperbolic} if $$\lambda\cdot l_X(H(\{0\}\times I))\leq\max\{l_X(H(\{m\}\times I)),l_X(H(\{-m\}\times I))\}.$$
		\item $X$ satisfies the {\bf hallways flare condition} if there exists $\lambda>1,m>0$ such that for every $b>0$ there exists $K(b)>0$ such that for any $b$-thin essential hallway of length $2m$ with girth $l_X(H(\{0\}\times I))\geq K(b)$ is $\lambda$-hyperbolic.
		\item An essential hallway of length 2m is {\em cone-bounded} if $H(i\times\partial I)$ lies in the cone locus, for $-m\leq i\leq m$.
		\item  $\widehat{X}$ satisfies the {\bf cone-bounded hallways strictly flare condition} if there exists $\lambda>1,m>0$ such that any cone-bounded essential hallway of length $2m$ is $\lambda$-hyperbolic.
	\end{enumerate}
\end{definition}
\begin{theorem}\cite[Theorem 4.6]{MjReeves}\label{mjreeves}
	Suppose $G$ admits a graph of relatively hyperbolic groups splitting that satisfies the following:\\
	$(1)$ the qi-embedded condition\\
	$(2)$ the strictly type-preserving condition\\
	$(3)$ the qi-preserving electrocution condition\\
	$(4)$ the induced tree of coned-off spaces satisfies the hallways flare condition\\
	$(5)$ the cone-bounded hallways strictly flare condition.\\
	Then $G$ is hyperbolic relative to a family $\mathcal{C}$ of maximal parabolic subgroups.
\end{theorem}
%\begin{theorem}\cite[Theorem 4.5]{MjReeves}\label{mjreeves}
%Let $X$ be a tree of relatively hyperbolic spaces satisfying\\
%$(1)$ the qi-embedded condition,\\
%$(2)$ the strictly type-preserving condition,\\
%$(3)$ the qi-preserving electrocution condition,\\
%$(4)$ the induced tree of coned-off spaces satisfies the hallways flare condition,\\
%$(5)$ the cone-bounded hallways strictly flare condition.\\
%Then $X$ is hyperbolic relative to the family $\CC$ of maximal cone-subtrees of horosphere-like spaces.
%\end{theorem}
\begin{definition}
	$\f:\T\to\T$ is {\bf $(\lambda,M)$-hyperbolic} if there exists $M>0,\lambda>1$, such that if $\rho$ is the fundamental segment of a hyperbolic element in $\T$ or a geodesic joining two non-free vertices, then $$(\ast\ast)\hspace{5mm}\lambda l_{\T}([\f^M(\rho)])\leq\max\{l_{\T}(\rho),l_{\T}([\f^{2M}(\rho)])\}.$$ 
\end{definition}
\noindent
This is equivalent to the flaring conditions in Theorem \ref{mjreeves} and we have noted this in the below lemma.
\begin{lemma}\label{equivalence-flare}
	$\mathcal{TC}(X)$ satisfies a hallways flare condition and the cone-bounded hallway strictly flare condition if and only if $\f:\T\to\T$ is hyperbolic.
\end{lemma}
\begin{proof}
	Suppose $\f:\T\to\T$ is $(\lambda,M)$-hyperbolic. Choose the smallest $k>0$ such that $\lambda^k>3$. Then $\f$ is $(3,2^kM)$-hyperbolic. Let $H:[-2^kM,2^kM]\times I\to\mathcal{TC}(X)$ be an essential, $b$-thin hallway. For each $i\in[-2^kM,2^kM]\cap\Z$, $\alpha_i:=H(\{i\}\times I)$ is a geodesic in a vertex space that we will denote by $X_{v_i}$. Recall that the vertex spaces are copies of $\T$. We denote the metric on $X_{v_i}$ by $d_i$. 
	
	Recall that the vertex spaces in $\mathcal{TC}(X)$ are uniformly properly embedded in it, so there exists an increasing function $\eta:[0,\infty)\to[0,\infty)$ such that for any $x,y\in\hat{X}_v$, if $d_{\mathcal{TC}}(x,y)\leq n$, then $d_{\X_v}(x,y)\leq\eta(n)$. Also, note that, for a vertex space $\X_v$ and some $x\in V(\X_v)$, $d_{\mathcal{TC}(X)}(x,\f(x))\leq 1$. Put $N=2^kM$. Let $l_0(\alpha_0)\geq7\eta((b+1)N)$. Let $\alpha_{-M}$ be a geodesic that joins free vertices. 
	By induction, we can show that $$l_0(\f^N(\alpha_{-N}))\geq l_0(\alpha_0)-2\eta(N(b+1)).$$
	\noindent
	Suppose $3l_0([\f^N(\alpha_{-N})])\leq l_{-N}(\alpha_{-N})$. Then,\begin{equation}
		%\begin{split}
		l_{-N}(\alpha_{-N})\geq3l_0(\alpha_0)-6\eta((b+1)N)       >2l_0(\alpha_0).
		%\end{split}
	\end{equation}
	Now suppose, $3l_0([\f^N(\alpha_{-N})])\leq l_{N}([\f^{2N}(\alpha_{-N})])$. Then as above, \begin{equation}l_{N}([\f^{2N}(\alpha_{-N})])>2l_0(\alpha_0). \end{equation}
	Now, let $\alpha_{-M}$ be a geodesic that joins a free vertex and a non-free vertex and suppose that this non-free vertex is stabilized by $H$. Then we denote by $\beta$ the concatenation of $\alpha_{-N}$ and $h\alpha_{-N}$ for some $h\in H$. Note that since $\f$ is a qi embedding and $\X_v$'s are uniformly hyperbolic, there exists some $k>0$ such that for $j=1,2$, $$2l_0(\f^{jN}(\alpha_{-N}))\geq l_0(\f^{jN}(\beta))\geq 2l_0(\f^{jN}(\alpha_{-N}))-2k.$$ 
	Also in this case, by induction, $l_0(\f^N(\alpha_{-N}))\geq l_0(\alpha_0)-\eta(N(b+1))$.\\
	Let $l_0(\alpha_0)>4\eta((b+1)N)+4k$. Now, suppose $3l_0([\f^N(\beta)])\leq l_{-N}(\beta)$. Then, \begin{equation}
		%\begin{split}
		2l_{-N}(\alpha_{-N})=l_{-N}(\beta)\geq 3l_0([\f^N(\beta)])\geq 6l_0(\f^N(\alpha_{-N}))-6k> 2l_0(\alpha_0).        
		%\end{split}
	\end{equation}
	\noindent
	Similarly, when $3l_0([\f^N(\beta)])\leq l_{N}([\f^{2N}(\beta)])$, we have \begin{equation}l_{N}([\f^{2N}(\alpha_{-N})])>2l_0(\alpha_0).\end{equation}
	Therefore, any $b$-thin essential hallway $H:[2^k,M,2^kM]\times I\to\mathcal{TC}(X)$ with $$K(b)=\min\{7\eta((b+1)N),4\eta((b+1)N)+4k\},$$ flares.
	%Suppose $H$ is a cone-bounded hallway. By the hyperbolicity of $\f$, since $\alpha_i$'s are geodesic segments joining non-free vertices, $$\lambda l_{X_{v_0}}(\alpha_0)\leq\max\{l_{X_{v_{-m}}}(\alpha_{-m}),l_{X_{v_m}}(\alpha_m)\}.$$
	\noindent
	Now, if $H$ is a cone-bounded hallway, then for all $i\in(-N,N]\cap\Z$, $\alpha_{i}=[\f^i(\alpha_{-N})]$, and thus, the cone-bounded hallways strictly flare condition is also satisfied.
	
	Conversely, it is easy to see that if $\mathcal{TC}(X)$ satisfies a hallways flare condition and cone-bounded hallways strictly flare condition, then $\f:\T\to\T$ is hyperbolic.\end{proof}

In the proof of the following result, the first part is the same as that of \cite[Proposition 2.18]{DahmaniLi}. 
\begin{theorem}\label{general}
	Let $G$ be a finitely generated group with a free factor system $\GG$, such that $G\ncong\Z$. Let $\phi:(G,\GG)\to(G,\GG)$ be an injective non-surjective endomorphism and $\PP$ be a $\phi$-invariant collection of subgroups of $G$ such that the following hold:\\
	(i) $\phi$ is fully irreducible relative to $\GG$.\\
	(ii) $\PP$ hyperbolically cones-off $\GG$.\\
	(iii) $\phi$ is strictly type preserving relative to $\PP$.\\
	%(iv) $\PP$ is transverse to legal paths of $\phi$.\\
	(iv) $\phi$ is atoroidal relative to $\PP$.\\
	(v) $\phi$ has no twinned subgroups for $\PP$.\\
	\noindent
	Suppose there exists $A>1$ such that for any legal path $\rho$ in $T$, we have $$l_{\T}([\hrho])\geq \frac{1}{A}l_T([\rho]).$$
	Then the ascending HNN extension $G\ast_{\phi}=\langle G,t\mid t^{-1}gt=\phi(g),g\in G\rangle$ is hyperbolic relative to a collection of maximal parabolic subgroups.
\end{theorem}
\begin{proof}
	If $\lambda_f=1$, it is easy to see that either $G=H$ with $\GG=\{[H]\}$ or $G=H\ast\ZZ$ with $\GG=\{[H\ast\ZZ]\}$. In these two cases, taking $\PP=\{[G]\}$, $G\ast_{\phi}$ is hyperbolic relative to itself.  So we assume that $\lambda_f>1$. By replacing $\lambda_f$ with a power, we assume that $\lambda_f>A$. By Lemma \ref{relhypG}, $G$ is hyperbolic relative to $\PP$. Let $\sigma$ either be a geodesic that is the fundamental segment of a hyperbolic element $g\in G$ in $\T$ or joins two non-free vertices in $\T$. For a geodesic $\rho\subset \T$, recall that $$\text{LEG}_{\T}(\rho)=\frac{\text{sum of legal segments of } \hrho \text{ of length at least } C(f)\text{ in }\T}{l_{\T}(\hrho)}.$$
	
	\noindent
	$(\ast)$ Suppose there exists $\epsilon>0,N_1>0$ such that $\text{LEG}_{\T}(\f^n(\sigma))\geq\epsilon$ for every $n\geq N_1$. Let $n\geq N_1$. Then, $$f^{n}(\sigma)=\sigma_1\cdot\sigma_2\cdots\sigma_m,$$ where $\sigma_i$ is a maximal legal segment, for $1\leq i\leq m$. Then by bounded cancellation lemma, we have, $$f^{2n}(\sigma)=\sigma'_1\cdots\sigma'_k,$$ where $f^n(\sigma_1)=\sigma'_1\cdot u_1$, $f^n(\sigma_i)=\overline{u_{i-1}}\cdot\sigma'_1\cdot u_i$, for $1<i<k$, and $f^n(\sigma_k)=\overline{u_{k-1}}\cdot\sigma'_k$ and $l_T(u_i)\leq C_{\ref{bcl}}$.
	
	\noindent
	Let $\Sigma=\{\sigma_i\mid l_{\T}(\hsigma_i)\geq C(f)\}$. Then by Lemma \ref{bcl-appl}, $$l_{\T}(\f^{2n}(\hsigma))\geq\sum_{\sigma_i\in\Sigma}l_{\T}(\hsigma'_i)\geq\nu\big(\frac{\lambda_f}{C_{\ref{transversality}}}\big)^{n}\sum_{\sigma_i\in\Sigma}l_T(\hsigma_i)\geq\nu\big(\frac{\lambda_f}{C_{\ref{transversality}}}\big)^{n}\left(\epsilon\cdot l_{\T}(\f^{n}(\hsigma))\right).$$
	\noindent
	Therefore, for every $n\geq N_1$, \begin{equation}\label{1}
		l_{\T}(\f^{2n}(\hsigma))\geq\nu\big(\frac{\lambda_f}{C_{\ref{transversality}}}\big)^{n}\left(\epsilon\cdot l_{\T}(\f^{n}(\hsigma))\right).    
	\end{equation}
	
	\noindent
	Now suppose $(\ast)$ does not hold. Then as $n\to\infty$, $\text{LEG}_{\T}(\f^n(\hsigma))\to0$. So, there exists $N_2>0$ such that for $n\geq N_2$ and any legal segment $\hat{\mu}$ in $\f^n(\hsigma)$, $l_{\T}(\hat{\mu})<C(f)$.\\
	Suppose $i(\f^n(\hsigma))<M_{\ref{nielsen}}$. Since $f$ is a train track map, the number of illegal turns does not increase under $f$, therefore, $i(\f^n(\hsigma))=i(\f^{2n}(\hsigma))$.\\
	We have, $$l_{\T}(\f^n(\hsigma))\leq i(\f^n(\hsigma))\cdot C(f)<M_{\ref{nielsen}}\cdot C(f).$$ For any $n\geq N_2$, we can take $\f^n(\hsigma)=\sigma_1\cdot\sigma_2\cdots\sigma_m$ where each $\sigma_i$ is a maximal legal subsegment and $m<M_{\ref{nielsen}}$. 
	\begin{comment}
	Then for every $k\geq 1$, 
	\begin{equation*}
	\begin{split}
	l_{\T}(\f^{(k+1)n}(\sigma))&\geq\sum_{i=1}^ml_{\T}(\f^{kn}(\sigma_i))-2(m-1)C_{\ref{bcl}}\\
	&\geq\lambda^{kn}\sum_{i=1}^ml_T(\sigma_i)-2(m-1)C_{\ref{bcl}}\\
	&>\lambda^{kn}\sum_{i=1}^ml_T(\sigma_i)-2M_{\ref{nielsen}}C_{\ref{bcl}}\\
	\end{split}
	\end{equation*}
	\end{comment}
	Then, \begin{equation*}
		\begin{split}
			l_{\T}(\f^{2n}(\sigma))&\geq\sum_{i=1}^ml_{\T}(\f^{n}(\sigma_i))-2(m-1)C_{\ref{bcl}}\\
			&\geq\lambda^{n}\sum_{i=1}^ml_T(\sigma_i)-2(m-1)C_{\ref{bcl}}\\
			&>\lambda^{n}\sum_{i=1}^ml_T(\sigma_i)-2M_{\ref{nielsen}}C_{\ref{bcl}}\\
		\end{split}
	\end{equation*}
	\noindent
	Choose $N_3\geq N_2$ large enough that for all $n\geq N_3$, $\lambda^n\sum_{i=1}^ml_T(\sigma_i)-2M_{\ref{nielsen}}C_{\ref{bcl}}>2M_{\ref{nielsen}}C_{\ref{bcl}}>2l_{\T}(\f^n(\hsigma))$.
	Thus, for $n\geq N_3$, \begin{equation}\label{2} 
		l_{\T}(\f^{2n}(\sigma))>2l_{\T}(\f^n(\hsigma)).
	\end{equation}
	\begin{comment} 
	\noindent
	Choose $k$ large enough that $\lambda^{kn}\sum_{i=1}^ml_T(\sigma_i)-2M_{\ref{nielsen}}C_{\ref{bcl}}>2M_{\ref{nielsen}}C_{\ref{bcl}}>2l_{\T}(\f^n(\hsigma))$.
	\end{comment}
	\noindent 
	Now, let $i(\f^n(\hsigma))\geq M_{\ref{nielsen}}$. Let $N_4\geq1$ such that for every $n\geq N_4$, $K_{\ref{length-illegal}}^2<\big(\frac{M_{\ref{nielsen}}+2}{M_{\ref{nielsen}}+1}\big)^{n}$. Then for every $n\geq N_4$ we have,
	\begin{equation}\label{3}
		l_{\T}(\hsigma)\geq K_{\ref{length-illegal}}^{-1}\big(\frac{M_{\ref{nielsen}}+2}{M_{\ref{nielsen}}+1}\big)^ni(\f^{n^2}(\hsigma))\geq K_{\ref{length-illegal}}^{-2}\big(\frac{M_{\ref{nielsen}}+2}{M_{\ref{nielsen}}+1}\big)^nl(\f^{n^2}(\hsigma)).   
	\end{equation}
	%Therefore, for every $n\geq\max\{N_2,N_3\}$,
	%$$l(\hsigma)\geq K_{\ref{length-illegal}}^{-2}\big(\frac{M_0+2}{M_0+1}\big)^nl(\f^{n^2}(\hsigma))$$.
	\noindent 
	Now, take $N=\max\{N_1,N_3,N_4\}$ and $$\lambda=\min\{\nu\epsilon\big(\frac{\lambda_f}{C_{\ref{transversality}}}\big)^{N^2},2,K_{\ref{length-illegal}}^{-2}\big(\frac{M_0+2}{M_0+1}\big)^N\}.$$
	$$\text{Then,\hspace{5mm}}\lambda\cdot l_{\T}(\f^{N^2}(\hsigma))\geq\max\{l_{\T}(\hsigma),l_{\T}(\f^{2N^2}(\hsigma))\}.$$ Thus, $\f$ is hyperbolic 
	%Let $F:[-N^2,N^2]\times I\to\TT\mathcal{C}(X)$ be given by $F(\{i\}\times I)=f^{i+N^2}(\sigma)$. This is an essential hallway of length $2N^2$. Taking $\sigma$ to be the fundamental segment of a hyperbolic element or $\sigma$ joins two non-free vertices, we have conditions $(4), (5)$ of Theorem \ref{mjreeves}. 
	and by Theorem \ref{mjreeves}, $G{\ast}_{\phi}$ is hyperbolic relative to a collection of maximal parabolic subgroups.    \end{proof}
\subsubsection{Maximal parabolic subgroups}\label{max-para} 
The maximal parabolic subgroups are the conjuagacy class of stabilizers the connected components in the cone-locus. We have $\PP=\{P_1,\ldots,P_q\}$.\\
{\bf Case 1:} Suppose for some $P_i\in\PP$, $\phi^j(P_i)\leq g_iP_ig_i^{-1}$ for some $j\geq1,g_i\in G\ast_{\phi}$. Then the component of the cone-locus containing $\nu_{P_i}$ is stabilized by the ascending HNN extension generated by $\{P_i,g_i^{-1}t^j\}$.\\
{\bf Case 2:} Suppose $\phi^j(P_i)\leq gP_kg^{-1}$ for some $P_k\in\PP,j\geq1,g\in G\ast_{\phi}$ such that $P_k$ is as in Case 1, i.e., $\phi^l(P_k)\leq g_kP_kg_k^{-1}$ for some $l\geq1, g_k\in G\ast_{\phi}$. Then the stabilizer of the component of the cone-locus containing $\nu_{P_i}$ is hyperbolic relative to the ascending HNN extension generated by $\{gP_kg^{-1},g_k^{-1}t^l\}$. 

Thus, there exists $\PP'\subseteq\PP$ and $K\geq1$ such that the maximal parabolic subgroups of $G\ast_{\phi^K}$ are ascending HNN extension of $P\in\PP'$.
\begin{theorem}\label{main}
	Let $G$ be a finitely generated group with a free factor system $\GG$, such that $G\ncong\Z$. Let $\phi:(G,\GG)\to(G,\GG)$ be an injective non-surjective endomorphism such that the following hold:\\
	(i) $\phi$ is fully irreducible relative to $\GG$.\\
	(ii) $\phi$ is atoroidal relative to $\GG$.\\
	(iii) $\phi$ is strictly type preserving relative to $\GG$.\\
	(iv) $\phi$ has no twinned subgroups.\\
	\noindent
	Then the ascending HNN extension $G\ast_{\phi}$ is hyperbolic relative to a collection of ascending HNN extensions of some elements in $\GG$.
\end{theorem}  
\noindent
This follows by taking $\PP=\GG$. Further, Theorem \ref{general} along with Proposition \ref{polynomial} and Lemma \ref{transversality} imply the following.
\begin{theorem}\label{main-poly}
	Let $G$ be a finitely generated group with a free factor system $\GG$, such that $G\ncong\Z$, and $\phi:(G,\GG)\to(G,\GG)$ be an injective endomorphism. Let $\PP$ be the collection of conjugacy classes of maximal polynomially growing subgroups for $\phi$ on a (any) tree in $\TT_{\GG}$. If $\phi$ is strictly type preserving relative to $\PP$, then the ascending HNN extension $G\ast_{\phi}=\langle G,t\mid t^{-1}gt=\phi(g),g\in G\rangle$ is hyperbolic relative to a collection of ascending HNN extension of some elements in $\PP$.
\end{theorem}
\begin{comment}
Recall that $p:\mathcal{TC}(X)\to\mathcal{T}$ is the tree of induced coned-off spaces with copies of $\widehat{T}=\mathcal{E}(T,\{T_P\}_{[P]\in\PP})$ as the vertex and edge spaces. We denote each vertex of the cone-locus by $(\nu(P),xG)$, where $[P]\in\PP$ and $\nu(P)$ is the cone point in the copy of $\T$ corresponding to $xG\in V(\mathcal{T})$. For $[P],[P']\in\PP, xG,yG\in V(\mathcal{T})$, there exists an edge from $(\nu(P),xG)$ to $(\nu(P'),yG)$ if there exists an edge joining $xG,yG$ in $\mathcal{T}$ and $f(T_P)\subset T_{P'}$. 

We have $\PP=\{[P_1],\ldots,[P_q]\}$. For $1\leq i\leq q$, let $Y_i$ be the connected component of the cone-locus containing the vertex $(\nu(P_i),G)$.
Let $G_i$ be the stabilizer of $Y_i$ in $G$. Note that any connected component of the cone-locus lies in the orbit of some $Y_i$ under the action of $G$ and thus, its stabilizer is a conjugate of $G_i$. Also, any conjugate of $G_i$ stabilizes a connected component of the cone-locus. So, the set maximal parabolic subgroups is given by $\{G_1,\ldots,G_q\}$.
\end{comment}

\subsection{Relative hyperbolicity of $\FF\ast_{\phi}$}
Let $\FF=\FF_n,n\geq 2$. Recall that an injective endomorphism $\phi:\FF\to\FF$ is {\it exponentially growing} if there exists $g\in\FF\setminus\{1\}$ and $C>1$ such that $||\phi^i(g)||\geq C||g||^i$ for all $i\geq 0$. 
\begin{lemma}\label{irred_ff}\cite[Lemma 1.4]{DahmaniLi}
	Let $G$ be a finitely presented group. Let $\phi:G\to G$ be an injective non-surjective endomorphism. Then there exists a proper $\phi$-invariant free factor system $\GG_0$ such that some power of $\phi$ is fully irreducible relative to $\GG_0$. 
\end{lemma}
\begin{proof}
	Suppose $\phi$ is not fully irreducible. Then there exists $i_1\geq1$ and $\GG_1\succ\{1\}$ such that $\phi^{i_1}(\GG_1)\preccurlyeq\GG_1$.    
	Proceeding inductively, if $\phi^{i_k}$ is not fully irreducible relative to $\GG_k$, there exists $i_{k+1}\geq i_k$ and $\GG_{k+1}\succ\GG_k$ such that $\phi^{i_{k+1}}(\GG_k)\preccurlyeq\GG_{k+1}$. Note that the Scott complexity of $\GG_{k+1}$ is strictly less than that of $\GG_k$. By Lemma \ref{scott_complexity}, there exists $N\in\mathbb{N}$ such that $\GG_N=\GG_{N+1}$. We take $\GG_0=\GG_N$. 
\end{proof}
\noindent
For an injective endomorphism $\phi$, we denote the free factor system obtained from Lemma \ref{irred_ff} by $\GG_{\phi}$. 
\begin{theorem}\label{main-free}
	Let $\phi:\FF\to\FF$ be an injective endomorphism that is exponentially growing. Let $\PP$ be the collection of conjugacy classes of maximal polynomially growing subgroups for $\phi$. If $\phi$ is strictly type preserving relative to $\PP$, then $\FF\ast_{\phi}=\langle\FF,t\mid t^{-1}gt=\phi(g),g\in\FF\rangle$ is hyperbolic relative to a collection of ascending HNN extension of some of the maximal polynomially growing subgroups.
\end{theorem}
\noindent
When $\phi$ is fully irreducible, by \cite[Theorem 6.3]{JPMnonsurj}, $\FF{\ast}_{\phi}$ is a hyperbolic group and $\PP=\{\emptyset\}$. If $\phi$ is not fully irreducible, then by Lemma \ref{irred_ff}, some power of $\phi$ is fully irreducible relative to a free factor system $\GG_0$. Then the result follows from Theorem \ref{main-poly}. And the maximal parabolic subgroups are as described in \ref{max-para}. 

\noindent
Mutanguha proved a combination theorem for relatively hyperbolicity.
\begin{theorem}\cite[Theorem 4.8]{JPMrelhyp}\label{jpmrelhyp}
	Suppose $\phi:\FF\to\FF$ is an injective, non-surjective endomorphism. Then, $\FF\ast_{\phi}$ is hyperbolic relative to a collection of free-by-cyclic groups if and only if $\FF\ast_{\phi}$ does not contain a subgroup isomorphic to $BS(1,d)$ for $d\geq 2$.
\end{theorem}
\noindent
{\em Sketch of the proof:} If $\phi$ is an expanding graph immersion, then by \cite[Theorem 6.3]{JPMhyperbolic}, $\FF\ast_{\phi}$ is hyperbolic if and only if $\FF\ast_{\phi}$ does not contain a subgroup isomorphic to $BS(1,d)$ for $d\geq 2$. If $\phi$ is not an expanding graph immersion, by \cite[Propositions 3.2.3, 3.2.4]{JPMrel} and \cite[Proposition 3.4.7]{JPMrel}, Mutanguha showed that there exists a $[\phi]$-fixed free factor system $\AA$ and a $[\phi]$-invariant free factor system $\AA^{\ast}$ such that $\AA\subset\AA^{\ast}$ and there exists some $k>0$ such that $\phi^k(\AA^{\ast})\preccurlyeq\AA$. It follows from \cite[Proposition 5.1.5]{JPMrel} and \cite[Proposition 5.2.3]{JPMrel} $\FF\ast_{\phi}$ is hyperbolic relative to a collection of maximal parabolic subgroups and these maximal parabolic subgroups are free-by-cyclic groups $F\rtimes\ZZ$, where each $F$ is a conjugate of some $A\in\AA$. \qed
	
	\section{Rapid decay property}\label{section 6}
	\noindent
We recall the following from \cite{chatterji-rdp}. Let $G$ be a finitely generated group with a finite generating set $S$.
%Recall that for a discrete group $G$, the group ring $$\C G=\{f:G\to\C\mid f(g)=0\text{ for all but finitely many }g\}$$ and $$l^2(G)=\{f:G\to\C\mid \sum_{g\in G}|f(g)|^2<\infty\}$$ is the Hilbert space with the inner product $\langle f,f'\rangle=\sum_{g\in G}f(g)\bar{f'(g)}$.
\noindent
The group $G$ acts on $l^2(G)$ by left regular representations. \begin{equation*} 
	\begin{split}
		& G\to \mathcal{U}(l^2(G))\\
		& g\mapsto \{\xi\mapsto g(\xi):=\delta_g\ast\xi\}
	\end{split}
\end{equation*}
\noindent
The left regular representation of $G$ induces an algebra homomorphism $\C G\to \BB(l^2(G))$.
%\begin{equation*} 
%\begin{split}
%& \\
%& f\mapsto \sum_{g\in G}f(g)\cdot\delta_g
%\end{split}
%\end{equation*}
%\noindent
For any $f\in\C G$, let $||f||_{\ast}$ denote the operator norm of the image of $f$ under this homomorphism. So, we have, $$||f||_{\ast}=\sup\{||f\ast\xi||_2\mid\xi\in l^2(G),||\xi||_2=1\}.$$
\begin{definition}
	Let $G$ be a finitely generated group with a finite generating set $S$. The group $G$ has the {\bf rapid decay property} if there exists constants $C,D$ such that for any $N\in\NN$, and any $f\in\C G$ that is supported on elements with word length less than $N$, $$||f||_{\ast}\leq C\cdot N^D||f||_2.$$
\end{definition}
\noindent
It follows from \cite{rdp-hyp,jolis-hyp} that hyperbolic groups have the rapid decay property. Further, relatively hyperbolic groups have the rapid decay property provided the maximal parabolic subgroups also have this property.
\begin{theorem}\cite[Theorem 1.1]{rdp-relhyp}\label{rdp-relhyp}
	Let $G$ be hyperbolic relative to a collection of finitely generated subgroups $\{H_1,\ldots,H_m\}$. Then $G$ has the rapid decay property if and only if each of the groups $H_1,\ldots,H_m$ has the rapid decay property.
\end{theorem}
\noindent
In \cite{gautero-lustig-rdp}, Gautero and Lustig showed the following.
\begin{theorem}\cite[Corollary 1.2]{gautero-lustig-rdp}\label{gautero-lustig-rdp}
	For any $\phi\in Aut(\FF_n)$, the group $\FF_n\rtimes_{\phi}\ZZ$ has the rapid decay property.
\end{theorem}
\noindent
We get an analogous result for ascending HNN extension of free groups.
\begin{theorem}\label{rapiddecay}
	Let $\phi:\FF\to\FF$ be an injective endomorphism. Then $\FF\ast_{\phi}$ does not have the rapid decay property if and only if $\FF\ast_{\phi}$ contains $BS(1,d)$ for some $d\geq 2$.
\end{theorem}
\begin{proof}
	Suppose $\FF\ast_{\phi}$ does not contain $BS(1,d)$ for any $d\geq 2$. Then by Theorem \ref{jpmrelhyp}, $\FF\ast_{\phi}$ is hyperbolic relative to a collection of free-by-cyclic subgroups and by Theorems \ref{gautero-lustig-rdp} and \ref{rdp-relhyp}, $\FF\ast_{\phi}$ has the rapid decay property.
	
	Now, suppose $\FF\ast_{\phi}$ contains $BS(1,d)$ for some $d\geq 2$. For $d\geq 2$, $BS(1,d)$ is an amenable group with exponential growth. Then, $BS(1,d)$ does not have the rapid decay property since an amenable group has the rapid decay property if and only if it has polynomial growth (\cite[Corollay 3.1.8]{jolis-hyp}). Also, a subgroup of a group with the rapid decay property also has this property (\cite[Proposition 2.1.1]{jolis-hyp}. Therefore, $\FF\ast_{\phi}$ does not have the rapid decay property.     
\end{proof}    
	
	\bibliography{Combi}
	\bibliographystyle{amsalpha}

\end{document}